\documentclass[12pt]{article}
\usepackage{geometry}
\usepackage{amsfonts,amsmath,latexsym,amssymb,amstext}
\usepackage[francais]{babel}
\usepackage{times}
\usepackage{eufrak}
\usepackage{synttree}
\usepackage{bussproofs}
\newdimen\proofrulebreadth \proofrulebreadth=.05em
\newdimen\proofdotseparation \proofdotseparation=1.25ex
\newdimen\proofrulebaseline \proofrulebaseline=2ex
\newcount\proofdotnumber \proofdotnumber=3
\let\then\relax
\def\hfi{\hskip0pt plus.0001fil}
\mathchardef\squigto="3A3B
\newif\ifinsideprooftree\insideprooftreefalse
\newif\ifonleftofproofrule\onleftofproofrulefalse
\newif\ifproofdots\proofdotsfalse
\newif\ifdoubleproof\doubleprooffalse
\let\wereinproofbit\relax
\newdimen\shortenproofleft
\newdimen\shortenproofright
\newdimen\proofbelowshift
\newbox\proofabove
\newbox\proofbelow
\newbox\proofrulename
\def\shiftproofbelow{\let\next\relax\afterassignment\setshiftproofbelow\dimen0 }
\def\shiftproofbelowneg{\def\next{\multiply\dimen0 by-1 }%
\afterassignment\setshiftproofbelow\dimen0 }
\def\setshiftproofbelow{\next\proofbelowshift=\dimen0 }
\def\setproofrulebreadth{\proofrulebreadth}

\def\prooftree{
\ifnum	\lastpenalty=1
\then	\unpenalty
\else	\onleftofproofrulefalse
\fi
\ifonleftofproofrule
\else	\ifinsideprooftree
	\then	\hskip.5em plus1fil
	\fi
\fi
\bgroup
\setbox\proofbelow=\hbox{}\setbox\proofrulename=\hbox{}%
\let\justifies\proofover\let\leadsto\proofoverdots\let\Justifies\proofoverdbl
\let\using\proofusing\let\[\prooftree
\ifinsideprooftree\let\]\endprooftree\fi
\proofdotsfalse\doubleprooffalse
\let\thickness\setproofrulebreadth
\let\shiftright\shiftproofbelow \let\shift\shiftproofbelow
\let\shiftleft\shiftproofbelowneg
\let\ifwasinsideprooftree\ifinsideprooftree
\insideprooftreetrue
\setbox\proofabove=\hbox\bgroup$\displaystyle 
\let\wereinproofbit\prooftree
\shortenproofleft=0pt \shortenproofright=0pt \proofbelowshift=0pt
\onleftofproofruletrue\penalty1
}

\def\eproofbit{
\ifx	\wereinproofbit\prooftree
\then	\ifcase	\lastpenalty
	\then	\shortenproofright=0pt	
	\or	\unpenalty\hfil		
	\or	\unpenalty\unskip	
	\else	\shortenproofright=0pt	
	\fi
\fi
\global\dimen0=\shortenproofleft
\global\dimen1=\shortenproofright
\global\dimen2=\proofrulebreadth
\global\dimen3=\proofbelowshift
\global\dimen4=\proofdotseparation
\global\count10=\proofdotnumber
$\egroup  
\shortenproofleft=\dimen0
\shortenproofright=\dimen1
\proofrulebreadth=\dimen2
\proofbelowshift=\dimen3
\proofdotseparation=\dimen4
\proofdotnumber=\count10
}

\def\proofover{
\eproofbit 
\setbox\proofbelow=\hbox\bgroup 
\let\wereinproofbit\proofover
$\displaystyle
}%
\def\proofoverdbl{
\eproofbit 
\doubleprooftrue
\setbox\proofbelow=\hbox\bgroup 
\let\wereinproofbit\proofoverdbl
$\displaystyle
}%
\def\proofoverdots{
\eproofbit 
\proofdotstrue
\setbox\proofbelow=\hbox\bgroup 
\let\wereinproofbit\proofoverdots
$\displaystyle
}%
\def\proofusing{
\eproofbit 
\setbox\proofrulename=\hbox\bgroup 
\let\wereinproofbit\proofusing
\kern0.3em$
}

\def\endprooftree{
\eproofbit 
  \dimen5 =0pt
\dimen0=\wd\proofabove \advance\dimen0-\shortenproofleft
\advance\dimen0-\shortenproofright
\dimen1=.5\dimen0 \advance\dimen1-.5\wd\proofbelow
\dimen4=\dimen1
\advance\dimen1\proofbelowshift \advance\dimen4-\proofbelowshift
\ifdim	\dimen1<0pt
\then	\advance\shortenproofleft\dimen1
	\advance\dimen0-\dimen1
	\dimen1=0pt
	\ifdim  \shortenproofleft<0pt
        \then   \setbox\proofabove=\hbox{%
			\kern-\shortenproofleft\unhbox\proofabove}%
                \shortenproofleft=0pt
        \fi
\fi
\ifdim	\dimen4<0pt
\then	\advance\shortenproofright\dimen4
	\advance\dimen0-\dimen4
	\dimen4=0pt
\fi
\ifdim	\shortenproofright<\wd\proofrulename
\then	\shortenproofright=\wd\proofrulename
\fi
\dimen2=\shortenproofleft \advance\dimen2 by\dimen1
\dimen3=\shortenproofright\advance\dimen3 by\dimen4
\ifproofdots
\then
	\dimen6=\shortenproofleft \advance\dimen6 .5\dimen0
	\setbox1=\vbox to\proofdotseparation{\vss\hbox{$\cdot$}\vss}%
	\setbox0=\hbox{%
		\advance\dimen6-.5\wd1
		\kern\dimen6
		$\vcenter to\proofdotnumber\proofdotseparation
			{\leaders\box1\vfill}$%
		\unhbox\proofrulename}%
\else	\dimen6=\fontdimen22\the\textfont2 
	\dimen7=\dimen6
	\advance\dimen6by.5\proofrulebreadth
	\advance\dimen7by-.5\proofrulebreadth
	\setbox0=\hbox{%
		\kern\shortenproofleft
		\ifdoubleproof
		\then	\hbox to\dimen0{%
			$\mathsurround0pt\mathord=\mkern-6mu%
			\cleaders\hbox{$\mkern-2mu=\mkern-2mu$}\hfill
			\mkern-6mu\mathord=$}%
		\else	\vrule height\dimen6 depth-\dimen7 width\dimen0
		\fi
		\unhbox\proofrulename}%
	\ht0=\dimen6 \dp0=-\dimen7
\fi
\let\doll\relax
\ifwasinsideprooftree
\then	\let\VBOX\vbox
\else	\ifmmode\else$\let\doll=$\fi
	\let\VBOX\vcenter
\fi
\VBOX	{\baselineskip\proofrulebaseline \lineskip.2ex
	\expandafter\lineskiplimit\ifproofdots0ex\else-0.6ex\fi
	\hbox	spread\dimen5	{\hfi\unhbox\proofabove\hfi}%
	\hbox{\box0}%
	\hbox	{\kern\dimen2 \box\proofbelow}}\doll%
\global\dimen2=\dimen2
\global\dimen3=\dimen3
\egroup 
\ifonleftofproofrule
\then	\shortenproofleft=\dimen2
\fi
\shortenproofright=\dimen3
\onleftofproofrulefalse
\ifinsideprooftree
\then	\hskip.5em plus 1fil \penalty2
\fi
}

\newcommand\seq\vdash 
   
\newcommand\ts\otimes
\newcommand\ra\rightarrow
\newcommand\Ra\Rightarrow
\newcommand\llts\otimes
\newcommand\pt\bullet 
\newcommand\lts\bullet 
\newcommand\lto{\mathbin{\backslash}}

\newcommand{\linefle}{\mathrel{-\hspace{-0.2em} \circ }}

\newcommand{\pop}{\mbox{\mbox{--$\circ$}}}
\newcommand\dw\downarrow
\newtheorem{definition}{D\'efinition}

\begin{document}
%
%
%
%
%
%
\newcounter{treecount}
\newcounter{branchcount}
\setcounter{treecount}{0}
\newsavebox{\parentbox}
\newsavebox{\treebox}
\newsavebox{\treeboxone}
\newsavebox{\treeboxtwo}
\newsavebox{\treeboxthree}
\newsavebox{\treeboxfour}
\newsavebox{\treeboxfive}
\newsavebox{\treeboxsix}
\newsavebox{\treeboxseven}
\newsavebox{\treeboxeight}
\newsavebox{\treeboxnine}
\newsavebox{\treeboxten}
\newsavebox{\treeboxeleven}
\newsavebox{\treeboxtwelve}
\newsavebox{\treeboxthirteen}
\newsavebox{\treeboxfourteen}
\newsavebox{\treeboxfifteen}
\newsavebox{\treeboxsixteen}
\newsavebox{\treeboxseventeen}
\newsavebox{\treeboxeighteen}
\newsavebox{\treeboxnineteen}
\newsavebox{\treeboxtwenty}
\newlength{\treeoffsetone}
\newlength{\treeoffsettwo}
\newlength{\treeoffsetthree}
\newlength{\treeoffsetfour}
\newlength{\treeoffsetfive}
\newlength{\treeoffsetsix}
\newlength{\treeoffsetseven}
\newlength{\treeoffseteight}
\newlength{\treeoffsetnine}
\newlength{\treeoffsetten}
\newlength{\treeoffseteleven}
\newlength{\treeoffsettwelve}
\newlength{\treeoffsetthirteen}
\newlength{\treeoffsetfourteen}
\newlength{\treeoffsetfifteen}
\newlength{\treeoffsetsixteen}
\newlength{\treeoffsetseventeen}
\newlength{\treeoffseteighteen}
\newlength{\treeoffsetnineteen}
\newlength{\treeoffsettwenty}

\newlength{\treeshiftone}
\newlength{\treeshifttwo}
\newlength{\treeshiftthree}
\newlength{\treeshiftfour}
\newlength{\treeshiftfive}
\newlength{\treeshiftsix}
\newlength{\treeshiftseven}
\newlength{\treeshifteight}
\newlength{\treeshiftnine}
\newlength{\treeshiftten}
\newlength{\treeshifteleven}
\newlength{\treeshifttwelve}
\newlength{\treeshiftthirteen}
\newlength{\treeshiftfourteen}
\newlength{\treeshiftfifteen}
\newlength{\treeshiftsixteen}
\newlength{\treeshiftseventeen}
\newlength{\treeshifteighteen}
\newlength{\treeshiftnineteen}
\newlength{\treeshifttwenty}
\newlength{\treewidthone}
\newlength{\treewidthtwo}
\newlength{\treewidththree}
\newlength{\treewidthfour}
\newlength{\treewidthfive}
\newlength{\treewidthsix}
\newlength{\treewidthseven}
\newlength{\treewidtheight}
\newlength{\treewidthnine}
\newlength{\treewidthten}
\newlength{\treewidtheleven}
\newlength{\treewidthtwelve}
\newlength{\treewidththirteen}
\newlength{\treewidthfourteen}
\newlength{\treewidthfifteen}
\newlength{\treewidthsixteen}
\newlength{\treewidthseventeen}
\newlength{\treewidtheighteen}
\newlength{\treewidthnineteen}
\newlength{\treewidthtwenty}
\newlength{\daughteroffsetone}
\newlength{\daughteroffsettwo}
\newlength{\daughteroffsetthree}
\newlength{\daughteroffsetfour}
\newlength{\branchwidthone}
\newlength{\branchwidthtwo}
\newlength{\branchwidththree}
\newlength{\branchwidthfour}
\newlength{\parentoffset}
\newlength{\treeoffset}
\newlength{\daughteroffset}
\newlength{\branchwidth}
\newlength{\parentwidth}
\newlength{\treewidth}
\newcommand{\ontop}[1]{\begin{tabular}{c}#1\end{tabular}}
\newcommand{\poptree}{%
\ifnum\value{treecount}=0\typeout{QobiTeX warning---Tree stack underflow}\fi%
\addtocounter{treecount}{-1}%
\setlength{\treeoffsettwo}{\treeoffsetthree}%
\setlength{\treeoffsetthree}{\treeoffsetfour}%
\setlength{\treeoffsetfour}{\treeoffsetfive}%
\setlength{\treeoffsetfive}{\treeoffsetsix}%
\setlength{\treeoffsetsix}{\treeoffsetseven}%
\setlength{\treeoffsetseven}{\treeoffseteight}%
\setlength{\treeoffseteight}{\treeoffsetnine}%
\setlength{\treeoffsetnine}{\treeoffsetten}%
\setlength{\treeoffsetten}{\treeoffseteleven}%
\setlength{\treeoffseteleven}{\treeoffsettwelve}%
\setlength{\treeoffsettwelve}{\treeoffsetthirteen}%
\setlength{\treeoffsetthirteen}{\treeoffsetfourteen}%
\setlength{\treeoffsetfourteen}{\treeoffsetfifteen}%
\setlength{\treeoffsetfifteen}{\treeoffsetsixteen}%
\setlength{\treeoffsetsixteen}{\treeoffsetseventeen}%
\setlength{\treeoffsetseventeen}{\treeoffseteighteen}%
\setlength{\treeoffseteighteen}{\treeoffsetnineteen}%
\setlength{\treeoffsetnineteen}{\treeoffsettwenty}%
\setlength{\treeshifttwo}{\treeshiftthree}%
\setlength{\treeshiftthree}{\treeshiftfour}%
\setlength{\treeshiftfour}{\treeshiftfive}%
\setlength{\treeshiftfive}{\treeshiftsix}%
\setlength{\treeshiftsix}{\treeshiftseven}%
\setlength{\treeshiftseven}{\treeshifteight}%
\setlength{\treeshifteight}{\treeshiftnine}%
\setlength{\treeshiftnine}{\treeshiftten}%
\setlength{\treeshiftten}{\treeshifteleven}%
\setlength{\treeshifteleven}{\treeshifttwelve}%
\setlength{\treeshifttwelve}{\treeshiftthirteen}%
\setlength{\treeshiftthirteen}{\treeshiftfourteen}%
\setlength{\treeshiftfourteen}{\treeshiftfifteen}%
\setlength{\treeshiftfifteen}{\treeshiftsixteen}%
\setlength{\treeshiftsixteen}{\treeshiftseventeen}%
\setlength{\treeshiftseventeen}{\treeshifteighteen}%
\setlength{\treeshifteighteen}{\treeshiftnineteen}%
\setlength{\treeshiftnineteen}{\treeshifttwenty}%
\setlength{\treewidthtwo}{\treewidththree}%
\setlength{\treewidththree}{\treewidthfour}%
\setlength{\treewidthfour}{\treewidthfive}%
\setlength{\treewidthfive}{\treewidthsix}%
\setlength{\treewidthsix}{\treewidthseven}%
\setlength{\treewidthseven}{\treewidtheight}%
\setlength{\treewidtheight}{\treewidthnine}%
\setlength{\treewidthnine}{\treewidthten}%
\setlength{\treewidthten}{\treewidtheleven}%
\setlength{\treewidtheleven}{\treewidthtwelve}%
\setlength{\treewidthtwelve}{\treewidththirteen}%
\setlength{\treewidththirteen}{\treewidthfourteen}%
\setlength{\treewidthfourteen}{\treewidthfifteen}%
\setlength{\treewidthfifteen}{\treewidthsixteen}%
\setlength{\treewidthsixteen}{\treewidthseventeen}%
\setlength{\treewidthseventeen}{\treewidtheighteen}%
\setlength{\treewidtheighteen}{\treewidthnineteen}%
\setlength{\treewidthnineteen}{\treewidthtwenty}%
\sbox{\treeboxtwo}{\usebox{\treeboxthree}}%
\sbox{\treeboxthree}{\usebox{\treeboxfour}}%
\sbox{\treeboxfour}{\usebox{\treeboxfive}}%
\sbox{\treeboxfive}{\usebox{\treeboxsix}}%
\sbox{\treeboxsix}{\usebox{\treeboxseven}}%
\sbox{\treeboxseven}{\usebox{\treeboxeight}}%
\sbox{\treeboxeight}{\usebox{\treeboxnine}}%
\sbox{\treeboxnine}{\usebox{\treeboxten}}%
\sbox{\treeboxten}{\usebox{\treeboxeleven}}%
\sbox{\treeboxeleven}{\usebox{\treeboxtwelve}}%
\sbox{\treeboxtwelve}{\usebox{\treeboxthirteen}}%
\sbox{\treeboxthirteen}{\usebox{\treeboxfourteen}}%
\sbox{\treeboxfourteen}{\usebox{\treeboxfifteen}}%
\sbox{\treeboxfifteen}{\usebox{\treeboxsixteen}}%
\sbox{\treeboxsixteen}{\usebox{\treeboxseventeen}}%
\sbox{\treeboxseventeen}{\usebox{\treeboxeighteen}}%
\sbox{\treeboxeighteen}{\usebox{\treeboxnineteen}}%
\sbox{\treeboxnineteen}{\usebox{\treeboxtwenty}}}
\newcommand{\leaf}[1]{%
\ifnum\value{treecount}=20\typeout{QobiTeX warning---Tree stack overflow}\fi%
\addtocounter{treecount}{1}%
\sbox{\treeboxtwenty}{\usebox{\treeboxnineteen}}%
\sbox{\treeboxnineteen}{\usebox{\treeboxeighteen}}%
\sbox{\treeboxeighteen}{\usebox{\treeboxseventeen}}%
\sbox{\treeboxseventeen}{\usebox{\treeboxsixteen}}%
\sbox{\treeboxsixteen}{\usebox{\treeboxfifteen}}%
\sbox{\treeboxfifteen}{\usebox{\treeboxfourteen}}%
\sbox{\treeboxfourteen}{\usebox{\treeboxthirteen}}%
\sbox{\treeboxthirteen}{\usebox{\treeboxtwelve}}%
\sbox{\treeboxtwelve}{\usebox{\treeboxeleven}}%
\sbox{\treeboxeleven}{\usebox{\treeboxten}}%
\sbox{\treeboxten}{\usebox{\treeboxnine}}%
\sbox{\treeboxnine}{\usebox{\treeboxeight}}%
\sbox{\treeboxeight}{\usebox{\treeboxseven}}%
\sbox{\treeboxseven}{\usebox{\treeboxsix}}%
\sbox{\treeboxsix}{\usebox{\treeboxfive}}%
\sbox{\treeboxfive}{\usebox{\treeboxfour}}%
\sbox{\treeboxfour}{\usebox{\treeboxthree}}%
\sbox{\treeboxthree}{\usebox{\treeboxtwo}}%
\sbox{\treeboxtwo}{\usebox{\treeboxone}}%
\sbox{\treeboxone}{\ontop{#1}}%
\sbox{\treeboxone}{\raisebox{-\ht\treeboxone}{\usebox{\treeboxone}}}%
\setlength{\treeoffsettwenty}{\treeoffsetnineteen}%
\setlength{\treeoffsetnineteen}{\treeoffseteighteen}%
\setlength{\treeoffseteighteen}{\treeoffsetseventeen}%
\setlength{\treeoffsetseventeen}{\treeoffsetsixteen}%
\setlength{\treeoffsetsixteen}{\treeoffsetfifteen}%
\setlength{\treeoffsetfifteen}{\treeoffsetfourteen}%
\setlength{\treeoffsetfourteen}{\treeoffsetthirteen}%
\setlength{\treeoffsetthirteen}{\treeoffsettwelve}%
\setlength{\treeoffsettwelve}{\treeoffseteleven}%
\setlength{\treeoffseteleven}{\treeoffsetten}%
\setlength{\treeoffsetten}{\treeoffsetnine}%
\setlength{\treeoffsetnine}{\treeoffseteight}%
\setlength{\treeoffseteight}{\treeoffsetseven}%
\setlength{\treeoffsetseven}{\treeoffsetsix}%
\setlength{\treeoffsetsix}{\treeoffsetfive}%
\setlength{\treeoffsetfive}{\treeoffsetfour}%
\setlength{\treeoffsetfour}{\treeoffsetthree}%
\setlength{\treeoffsetthree}{\treeoffsettwo}%
\setlength{\treeoffsettwo}{\treeoffsetone}%
\setlength{\treeoffsetone}{0.5\wd\treeboxone}%
\setlength{\treeshifttwenty}{\treeshiftnineteen}%
\setlength{\treeshiftnineteen}{\treeshifteighteen}%
\setlength{\treeshifteighteen}{\treeshiftseventeen}%
\setlength{\treeshiftseventeen}{\treeshiftsixteen}%
\setlength{\treeshiftsixteen}{\treeshiftfifteen}%
\setlength{\treeshiftfifteen}{\treeshiftfourteen}%
\setlength{\treeshiftfourteen}{\treeshiftthirteen}%
\setlength{\treeshiftthirteen}{\treeshifttwelve}%
\setlength{\treeshifttwelve}{\treeshifteleven}%
\setlength{\treeshifteleven}{\treeshiftten}%
\setlength{\treeshiftten}{\treeshiftnine}%
\setlength{\treeshiftnine}{\treeshifteight}%
\setlength{\treeshifteight}{\treeshiftseven}%
\setlength{\treeshiftseven}{\treeshiftsix}%
\setlength{\treeshiftsix}{\treeshiftfive}%
\setlength{\treeshiftfive}{\treeshiftfour}%
\setlength{\treeshiftfour}{\treeshiftthree}%
\setlength{\treeshiftthree}{\treeshifttwo}%
\setlength{\treeshifttwo}{\treeshiftone}%
\setlength{\treeshiftone}{0pt}%
\setlength{\treewidthtwenty}{\treewidthnineteen}%
\setlength{\treewidthnineteen}{\treewidtheighteen}%
\setlength{\treewidtheighteen}{\treewidthseventeen}%
\setlength{\treewidthseventeen}{\treewidthsixteen}%
\setlength{\treewidthsixteen}{\treewidthfifteen}%
\setlength{\treewidthfifteen}{\treewidthfourteen}%
\setlength{\treewidthfourteen}{\treewidththirteen}%
\setlength{\treewidththirteen}{\treewidthtwelve}%
\setlength{\treewidthtwelve}{\treewidtheleven}%
\setlength{\treewidtheleven}{\treewidthten}%
\setlength{\treewidthten}{\treewidthnine}%
\setlength{\treewidthnine}{\treewidtheight}%
\setlength{\treewidtheight}{\treewidthseven}%
\setlength{\treewidthseven}{\treewidthsix}%
\setlength{\treewidthsix}{\treewidthfive}%
\setlength{\treewidthfive}{\treewidthfour}%
\setlength{\treewidthfour}{\treewidththree}%
\setlength{\treewidththree}{\treewidthtwo}%
\setlength{\treewidthtwo}{\treewidthone}%
\setlength{\treewidthone}{\wd\treeboxone}}
\newcommand{\branch}[2]{%
\setcounter{branchcount}{#1}%
\ifnum\value{branchcount}=1\sbox{\parentbox}{\ontop{#2}}%
\setlength{\parentoffset}{\treeoffsetone}%
\addtolength{\parentoffset}{-0.5\wd\parentbox}%
\setlength{\daughteroffset}{0in}%
\ifdim\parentoffset<0in%
\setlength{\daughteroffset}{-\parentoffset}%
\setlength{\parentoffset}{0in}\fi%
\setlength{\parentwidth}{\parentoffset}%
\addtolength{\parentwidth}{\wd\parentbox}%
\setlength{\treeoffset}{\daughteroffset}%
\addtolength{\treeoffset}{\treeoffsetone}%
\setlength{\treewidth}{\wd\treeboxone}%
\addtolength{\treewidth}{\daughteroffset}%
\ifdim\treewidth<\parentwidth\setlength{\treewidth}{\parentwidth}\fi%
\sbox{\treebox}{\begin{minipage}{\treewidth}%
\begin{flushleft}%
\hspace*{\parentoffset}\usebox{\parentbox}\\
{\setlength{\unitlength}{2ex}%
\hspace*{\treeoffset}\begin{picture}(0,1)%
\put(0,0){\line(0,1){1}}%
\end{picture}}\\
\vspace{-\baselineskip}
\hspace*{\daughteroffset}%
\raisebox{-\ht\treeboxone}{\usebox{\treeboxone}}%
\end{flushleft}%
\end{minipage}}%
\setlength{\treeoffsetone}{\parentoffset}%
\addtolength{\treeoffsetone}{0.5\wd\parentbox}%
\setlength{\treeshiftone}{0pt}%
\setlength{\treewidthone}{\treewidth}%
\sbox{\treeboxone}{\usebox{\treebox}}%
\else\ifnum\value{branchcount}=2\sbox{\parentbox}{\ontop{#2}}%
\setlength{\branchwidthone}{\treewidthtwo}%
\addtolength{\branchwidthone}{\treeoffsetone}%
\addtolength{\branchwidthone}{-\treeshiftone}%
\addtolength{\branchwidthone}{-\treeoffsettwo}%
\setlength{\branchwidth}{\branchwidthone}%
\setlength{\daughteroffsetone}{\branchwidth}%
\addtolength{\daughteroffsetone}{-\branchwidthone}%
\addtolength{\daughteroffsetone}{-\treeshiftone}%
\setlength{\parentoffset}{-0.5\wd\parentbox}%
\addtolength{\parentoffset}{\treeoffsettwo}%
\addtolength{\parentoffset}{0.5\branchwidth}%
\setlength{\daughteroffset}{0in}%
\ifdim\parentoffset<0in%
\setlength{\daughteroffset}{-\parentoffset}%
\setlength{\parentoffset}{0in}\fi%
\setlength{\parentwidth}{\parentoffset}%
\addtolength{\parentwidth}{\wd\parentbox}%
\setlength{\treeoffset}{\daughteroffset}%
\addtolength{\treeoffset}{\treeoffsettwo}%
\setlength{\treewidth}{\wd\treeboxone}%
\addtolength{\treewidth}{\daughteroffsetone}%
\addtolength{\treewidth}{\treewidthtwo}%
\addtolength{\treewidth}{\daughteroffset}%
\ifdim\treewidth<\parentwidth\setlength{\treewidth}{\parentwidth}\fi%
\sbox{\treebox}{\begin{minipage}{\treewidth}%
\begin{flushleft}%
\hspace*{\parentoffset}\usebox{\parentbox}\\
{\setlength{\unitlength}{0.5\branchwidth}%
\hspace*{\treeoffset}\begin{picture}(2,0.5)%
\put(0,0){\line(2,1){1}}%
\put(2,0){\line(-2,1){1}}%
\end{picture}}\\
\vspace{-\baselineskip}
\hspace*{\daughteroffset}%
\makebox[\treewidthtwo][l]%
{\raisebox{-\ht\treeboxtwo}{\usebox{\treeboxtwo}}}%
\hspace*{\daughteroffsetone}%
\raisebox{-\ht\treeboxone}{\usebox{\treeboxone}}%
\end{flushleft}%
\end{minipage}}%
\setlength{\treeoffsetone}{\parentoffset}%
\addtolength{\treeoffsetone}{0.5\wd\parentbox}%
\setlength{\treeshiftone}{0pt}%
\setlength{\treewidthone}{\treewidth}%
\sbox{\treeboxone}{\usebox{\treebox}}\poptree%
\else\ifnum\value{branchcount}=3\sbox{\parentbox}{\ontop{#2}}%
\setlength{\branchwidthone}{\treewidthtwo}%
\addtolength{\branchwidthone}{\treeoffsetone}%
\addtolength{\branchwidthone}{-\treeshiftone}%
\addtolength{\branchwidthone}{-\treeoffsettwo}%
\setlength{\branchwidthtwo}{\treewidththree}%
\addtolength{\branchwidthtwo}{\treeoffsettwo}%
\addtolength{\branchwidthtwo}{-\treeshifttwo}%
\addtolength{\branchwidthtwo}{-\treeoffsetthree}%
\setlength{\branchwidth}{\branchwidthone}%
\ifdim\branchwidthtwo>\branchwidth%
\setlength{\branchwidth}{\branchwidthtwo}\fi%
\setlength{\daughteroffsetone}{\branchwidth}%
\addtolength{\daughteroffsetone}{-\branchwidthone}%
\addtolength{\daughteroffsetone}{-\treeshiftone}%
\setlength{\daughteroffsettwo}{\branchwidth}%
\addtolength{\daughteroffsettwo}{-\branchwidthtwo}%
\addtolength{\daughteroffsettwo}{-\treeshifttwo}%
\setlength{\parentoffset}{-0.5\wd\parentbox}%
\addtolength{\parentoffset}{\treeoffsetthree}%
\addtolength{\parentoffset}{\branchwidth}%
\setlength{\daughteroffset}{0in}%
\ifdim\parentoffset<0in%
\setlength{\daughteroffset}{-\parentoffset}%
\setlength{\parentoffset}{0in}\fi%
\setlength{\parentwidth}{\parentoffset}%
\addtolength{\parentwidth}{\wd\parentbox}%
\setlength{\treeoffset}{\daughteroffset}%
\addtolength{\treeoffset}{\treeoffsetthree}%
\setlength{\treewidth}{\wd\treeboxone}%
\addtolength{\treewidth}{\daughteroffsetone}%
\addtolength{\treewidth}{\treewidthtwo}%
\addtolength{\treewidth}{\daughteroffsettwo}%
\addtolength{\treewidth}{\treewidththree}%
\addtolength{\treewidth}{\daughteroffset}%
\ifdim\treewidth<\parentwidth\setlength{\treewidth}{\parentwidth}\fi%
\sbox{\treebox}{\begin{minipage}{\treewidth}%
\begin{flushleft}%
\hspace*{\parentoffset}\usebox{\parentbox}\\
{\setlength{\unitlength}{0.5\branchwidth}%
\hspace*{\treeoffset}\begin{picture}(4,1)%
\put(0,0){\line(2,1){2}}%
\put(2,0){\line(0,1){1}}%
\put(4,0){\line(-2,1){2}}%
\end{picture}}\\
\vspace{-\baselineskip}
\hspace*{\daughteroffset}%
\makebox[\treewidththree][l]%
{\raisebox{-\ht\treeboxthree}{\usebox{\treeboxthree}}}%
\hspace*{\daughteroffsettwo}%
\makebox[\treewidthtwo][l]%
{\raisebox{-\ht\treeboxtwo}{\usebox{\treeboxtwo}}}%
\hspace*{\daughteroffsetone}%
\raisebox{-\ht\treeboxone}{\usebox{\treeboxone}}%
\end{flushleft}%
\end{minipage}}%
\setlength{\treeoffsetone}{\parentoffset}%
\addtolength{\treeoffsetone}{0.5\wd\parentbox}%
\setlength{\treeshiftone}{0pt}%
\setlength{\treewidthone}{\treewidth}%
\sbox{\treeboxone}{\usebox{\treebox}}\poptree\poptree%
\else\ifnum\value{branchcount}=4\sbox{\parentbox}{\ontop{#2}}%
\setlength{\branchwidthone}{\treewidthtwo}%
\addtolength{\branchwidthone}{\treeoffsetone}%
\addtolength{\branchwidthone}{-\treeshiftone}%
\addtolength{\branchwidthone}{-\treeoffsettwo}%
\setlength{\branchwidthtwo}{\treewidththree}%
\addtolength{\branchwidthtwo}{\treeoffsettwo}%
\addtolength{\branchwidthtwo}{-\treeshifttwo}%
\addtolength{\branchwidthtwo}{-\treeoffsetthree}%
\setlength{\branchwidththree}{\treewidthfour}%
\addtolength{\branchwidththree}{\treeoffsetthree}%
\addtolength{\branchwidththree}{-\treeshiftthree}%
\addtolength{\branchwidththree}{-\treeoffsetfour}%
\setlength{\branchwidth}{\branchwidthone}%
\ifdim\branchwidthtwo>\branchwidth%
\setlength{\branchwidth}{\branchwidthtwo}\fi%
\ifdim\branchwidththree>\branchwidth%
\setlength{\branchwidth}{\branchwidththree}\fi%
\setlength{\daughteroffsetone}{\branchwidth}%
\addtolength{\daughteroffsetone}{-\branchwidthone}%
\addtolength{\daughteroffsetone}{-\treeshiftone}%
\setlength{\daughteroffsettwo}{\branchwidth}%
\addtolength{\daughteroffsettwo}{-\branchwidthtwo}%
\addtolength{\daughteroffsettwo}{-\treeshifttwo}%
\setlength{\daughteroffsetthree}{\branchwidth}%
\addtolength{\daughteroffsetthree}{-\branchwidththree}%
\addtolength{\daughteroffsetthree}{-\treeshiftthree}%
\setlength{\parentoffset}{-0.5\wd\parentbox}%
\addtolength{\parentoffset}{\treeoffsetfour}%
\addtolength{\parentoffset}{1.5\branchwidth}%
\setlength{\daughteroffset}{0in}%
\ifdim\parentoffset<0in%
\setlength{\daughteroffset}{-\parentoffset}%
\setlength{\parentoffset}{0in}\fi%
\setlength{\parentwidth}{\parentoffset}%
\addtolength{\parentwidth}{\wd\parentbox}%
\setlength{\treeoffset}{\daughteroffset}%
\addtolength{\treeoffset}{\treeoffsetfour}%
\setlength{\treewidth}{\wd\treeboxone}%
\addtolength{\treewidth}{\daughteroffsetone}%
\addtolength{\treewidth}{\treewidthtwo}%
\addtolength{\treewidth}{\daughteroffsettwo}%
\addtolength{\treewidth}{\treewidththree}%
\addtolength{\treewidth}{\daughteroffsetthree}%
\addtolength{\treewidth}{\treewidthfour}%
\addtolength{\treewidth}{\daughteroffset}%
\ifdim\treewidth<\parentwidth\setlength{\treewidth}{\parentwidth}\fi%
\sbox{\treebox}{\begin{minipage}{\treewidth}%
\begin{flushleft}%
\hspace*{\parentoffset}\usebox{\parentbox}\\
{\setlength{\unitlength}{0.5\branchwidth}%
\hspace*{\treeoffset}\begin{picture}(6,1)%
\put(0,0){\line(3,1){3}}%
\put(2,0){\line(1,1){1}}%
\put(4,0){\line(-1,1){1}}%
\put(6,0){\line(-3,1){3}}%
\end{picture}}\\
\vspace{-\baselineskip}
\hspace*{\daughteroffset}%
\makebox[\treewidthfour][l]%
{\raisebox{-\ht\treeboxfour}{\usebox{\treeboxfour}}}%
\hspace*{\daughteroffsetthree}%
\makebox[\treewidththree][l]%
{\raisebox{-\ht\treeboxthree}{\usebox{\treeboxthree}}}%
\hspace*{\daughteroffsettwo}%
\makebox[\treewidthtwo][l]%
{\raisebox{-\ht\treeboxtwo}{\usebox{\treeboxtwo}}}%
\hspace*{\daughteroffsetone}%
\raisebox{-\ht\treeboxone}{\usebox{\treeboxone}}%
\end{flushleft}%
\end{minipage}}%
\setlength{\treeoffsetone}{\parentoffset}%
\addtolength{\treeoffsetone}{0.5\wd\parentbox}%
\setlength{\treeshiftone}{0pt}%
\setlength{\treewidthone}{\treewidth}%
\sbox{\treeboxone}{\usebox{\treebox}}\poptree\poptree\poptree%
\else\ifnum\value{branchcount}=5\sbox{\parentbox}{\ontop{#2}}%
\setlength{\branchwidthone}{\treewidthtwo}%
\addtolength{\branchwidthone}{\treeoffsetone}%
\addtolength{\branchwidthone}{-\treeshiftone}%
\addtolength{\branchwidthone}{-\treeoffsettwo}%
\setlength{\branchwidthtwo}{\treewidththree}%
\addtolength{\branchwidthtwo}{\treeoffsettwo}%
\addtolength{\branchwidthtwo}{-\treeshifttwo}%
\addtolength{\branchwidthtwo}{-\treeoffsetthree}%
\setlength{\branchwidththree}{\treewidthfour}%
\addtolength{\branchwidththree}{\treeoffsetthree}%
\addtolength{\branchwidththree}{-\treeshiftthree}%
\addtolength{\branchwidththree}{-\treeoffsetfour}%
\setlength{\branchwidthfour}{\treewidthfive}%
\addtolength{\branchwidthfour}{\treeoffsetfour}%
\addtolength{\branchwidthfour}{-\treeshiftfour}%
\addtolength{\branchwidthfour}{-\treeoffsetfive}%
\setlength{\branchwidth}{\branchwidthone}%
\ifdim\branchwidthtwo>\branchwidth%
\setlength{\branchwidth}{\branchwidthtwo}\fi%
\ifdim\branchwidththree>\branchwidth%
\setlength{\branchwidth}{\branchwidththree}\fi%
\ifdim\branchwidthfour>\branchwidth%
\setlength{\branchwidth}{\branchwidthfour}\fi%
\setlength{\daughteroffsetone}{\branchwidth}%
\addtolength{\daughteroffsetone}{-\branchwidthone}%
\addtolength{\daughteroffsetone}{-\treeshiftone}%
\setlength{\daughteroffsettwo}{\branchwidth}%
\addtolength{\daughteroffsettwo}{-\branchwidthtwo}%
\addtolength{\daughteroffsettwo}{-\treeshifttwo}%
\setlength{\daughteroffsetthree}{\branchwidth}%
\addtolength{\daughteroffsetthree}{-\branchwidththree}%
\addtolength{\daughteroffsetthree}{-\treeshiftthree}%
\setlength{\daughteroffsetfour}{\branchwidth}%
\addtolength{\daughteroffsetfour}{-\branchwidthfour}%
\addtolength{\daughteroffsetfour}{-\treeshiftfour}%
\setlength{\parentoffset}{-0.5\wd\parentbox}%
\addtolength{\parentoffset}{\treeoffsetfive}%
\addtolength{\parentoffset}{2\branchwidth}%
\setlength{\daughteroffset}{0in}%
\ifdim\parentoffset<0in%
\setlength{\daughteroffset}{-\parentoffset}%
\setlength{\parentoffset}{0in}\fi%
\setlength{\parentwidth}{\parentoffset}%
\addtolength{\parentwidth}{\wd\parentbox}%
\setlength{\treeoffset}{\daughteroffset}%
\addtolength{\treeoffset}{\treeoffsetfive}%
\setlength{\treewidth}{\wd\treeboxone}%
\addtolength{\treewidth}{\daughteroffsetone}%
\addtolength{\treewidth}{\treewidthtwo}%
\addtolength{\treewidth}{\daughteroffsettwo}%
\addtolength{\treewidth}{\treewidththree}%
\addtolength{\treewidth}{\daughteroffsetthree}%
\addtolength{\treewidth}{\treewidthfour}%
\addtolength{\treewidth}{\daughteroffsetfour}%
\addtolength{\treewidth}{\treewidthfive}%
\addtolength{\treewidth}{\daughteroffset}%
\ifdim\treewidth<\parentwidth\setlength{\treewidth}{\parentwidth}\fi%
\sbox{\treebox}{\begin{minipage}{\treewidth}%
\begin{flushleft}%
\hspace*{\parentoffset}\usebox{\parentbox}\\
{\setlength{\unitlength}{0.5\branchwidth}%
\hspace*{\treeoffset}\begin{picture}(8,1)%
\put(0,0){\line(4,1){4}}%
\put(2,0){\line(2,1){2}}%
\put(4,0){\line(0,1){1}}%
\put(6,0){\line(-2,1){2}}%
\put(8,0){\line(-4,1){4}}%
\end{picture}}\\
\vspace{-\baselineskip}
\hspace*{\daughteroffset}%
\makebox[\treewidthfive][l]%
{\raisebox{-\ht\treeboxfour}{\usebox{\treeboxfive}}}%
\hspace*{\daughteroffsetfour}%
\makebox[\treewidthfour][l]%
{\raisebox{-\ht\treeboxfour}{\usebox{\treeboxfour}}}%
\hspace*{\daughteroffsetthree}%
\makebox[\treewidththree][l]%
{\raisebox{-\ht\treeboxthree}{\usebox{\treeboxthree}}}%
\hspace*{\daughteroffsettwo}%
\makebox[\treewidthtwo][l]%
{\raisebox{-\ht\treeboxtwo}{\usebox{\treeboxtwo}}}%
\hspace*{\daughteroffsetone}%
\raisebox{-\ht\treeboxone}{\usebox{\treeboxone}}%
\end{flushleft}%
\end{minipage}}%
\setlength{\treeoffsetone}{\parentoffset}%
\addtolength{\treeoffsetone}{0.5\wd\parentbox}%
\setlength{\treeshiftone}{0pt}%
\setlength{\treewidthone}{\treewidth}%
\sbox{\treeboxone}{\usebox{\treebox}}\poptree\poptree\poptree\poptree%
\else\typeout{QobiTeX warning--- Can't handle #1 branching}\fi\fi\fi\fi\fi}
\newcommand{\faketreewidth}[1]{%
\sbox{\parentbox}{\ontop{#1}}%
\setlength{\treewidthone}{0.5\wd\parentbox}%
\addtolength{\treewidthone}{\treeoffsetone}%
\setlength{\treeshiftone}{\treeoffsetone}%
\addtolength{\treeshiftone}{-0.5\wd\parentbox}}
\newcommand{\tree}{%
\usebox{\treeboxone}
\setlength{\treeoffsetone}{\treeoffsettwo}%
\sbox{\treeboxone}{\usebox{\treeboxtwo}}%
\poptree}

\title{Ludique: une logique sans axiome d'identit\'e?}
\author{Alain Lecomte\footnote{UMR "Structures Formelles de la Langue", CNRS-Universit\'e Paris 8 - Vincennes-Saint-Denis}  ~\footnote{en collaboration avec Myriam Quatrini, Institut de Math\'ematiques de Luminy}}
\date{}
\maketitle
\section{Introduction : l'importance des r\`egles structurelles}
\subsection{Rappels de logique lin\'eaire}
Depuis qu'on s'int\'eresse aux logiques sous-structurelles (\cite{SchDos93}), on \'etudie les effets des variations dans la pr\'esence et l'absence des r\`egles structurelles de la logique, sur le syst\`eme obtenu. Curieusement, plus on \'elimine de r\`egles structurelles, plus on tend \`a s'\'eloigner de l'id\'ealisme des formules pour se rapprocher du concret mat\'eriel des ressources ou des inscriptions locales. On peut faire remonter \`a 1958 et \`a l'article particuli\`erement f\'econd de J. Lambek (\cite{Lam58}) le premier syst\`eme logique se privant des r\`egles de contraction, d'affaiblissement et de permutation dans le cadre d'une pr\'esentation de la logique en termes de calcul des s\'equents. La suppression des r\`egles de contraction et d'affaiblissement se trouve \^etre le ressort particulier qui permet \`a la logique lin\'eaire d'exister (\cite{Gi87, Gi95}). On notera toutefois que cette derni\`ere n'est pas seulement une logique "sous-structurelle" mais, parce que son but est de permettre l'analyse fine des preuves et des programmes en termes de s\'emantique d\'enotationnelle, elle est surtout un cadre \`a l'int\'erieur duquel on peut exprimer aussi bien la logique intuitionniste que la logique classique. On utilise pour ce faire les fameuses exponentielles ({\bf !, ?}) et on obtient par exemple une d\'ecomposition de l'implication classique, en: 
$$
A\Ra B~ \equiv~ {\bf !} A\linefle B
$$
\'equivalence qui signifie que pour passer du caract\`ere concret d'une ressource \`a l'id\'eation d'une formule, encore faut-il admettre un op\'erateur qui rende la ressource p\'erenne (r\^ole de "{\bf !}"). L'implication lin\'eaire $\linefle$ consomme bien alors un $A$, mais il en reste pour une r\'eutilisation ult\'erieure.
\subsection{Du traitement des paradoxes}
Jouer ainsi sur le passage du concret \`a l'abstrait permet de r\'eexaminer de vieilles questions et de les \'eclairer sous un jour nouveau. Ainsi,
la d\'ecomposition ci-dessus fournit un instrument d'analyse int\'eressant pour l'\'etude des {\it paradoxes}. Certains auteurs (Shirohata, 1993) ont en effet mis en \'evidence la dissolution du paradoxe de Russell dans le cadre lin\'eaire. Soit $E = \{X; X\not\in X\}$, de $E\in E$ on d\'eduit $E\not\in E$, mais alors la pr\'emisse est consomm\'ee: elle n'est plus l\`a pour impliquer une contradiction. On est simplement conduit \`a envisager une suite infinie et alternative d'instants o\`u on asserte soit l'un soit l'autre de $E\in E$ et de $E\not\in E$. M\^eme raisonnement pour le paradoxe du Menteur. La contradiction ne r\'eappara{\^\i}t que pour les formules affect\'ees du symbole de r\'eutilisation d'une ressource: {\bf !}. L'analyse lin\'eaire nous permet donc d'op\'erer au sein des connaissances, une distinction entre celles qui sont {\it transitoires} et celles qui sont {\it p\'erennes}. Seules ces derni\`eres sont susceptibles d'entra{\^\i}ner des contradictions.\\
Le genre de syst\`eme obtenu par ce jeu sur les r\`egles structurelles est justement qualifi\'e de {\it sensible aux ressources}, nous devons  comprendre par l\`a qu'il est sensible \`a la {\it quantit\'e} voire \`a {\it l'ordre} dans lequel sont donn\'ees les ressources disponibles.
\subsection{Les logiques sensibles aux ressources et la langue}
La suppression de la r\`egle de permutation est \`a l'origine d'une sensibilit\'e \`a {\it l'ordre}: \`a condition de se cantonner dans la partie multiplicative du calcul et de restreindre celui-ci \`a des s\'equents intuitionistes, on obtient un syst\`eme lin\'eaire non commutatif correspondant au calcul de Lambek (avec seulement la diff\'erence que, dans ce dernier, on ajoute la contrainte d'ant\'ec\'edent non vide afin d'\'eviter la pr\'esence de types $A/A$ (ou $A\lto A$) qui ne seraient associ\'es \`a aucun \'el\'ement mat\'eriel). Cette sensibilit\'e \`a l'ordre est \`a la source de l'emploi de ce calcul en {\it linguistique}. La {\it syntaxe des langues naturelles} est \`a premi\`ere vue en effet l'exemple le plus typique d'un syst\`eme sensible \`a la fois \`a la quantit\'e et \`a l'ordre des ressources. Au niveau le plus \'el\'ementaire, un sujet ou un objet ne se consomment qu'une fois et on les d\'etermine \`a partir des configurations ordonn\'ees o\`u ils figurent dans l'analyse de la phrase\footnote{Du moins dans les langues dites configurationnelles}.\\
 Le programme minimaliste de Chomsky contient ainsi l'id\'ee que l'une des op\'erations fondamentales de la syntaxe, {\it Move} dite aussi depuis peu {\it Internal Merge}, doit son existence \`a la n\'ecessit\'e de v\'erifier la pr\'esence de traits formels (non interpr\'etables) sur des constituants ou des t\^etes lexicales. Lorsqu'on doit v\'erifier qu'un constituant poss\`ede bien le trait requis, on le d\'eplace en une certaine position o\`u peuvent s'\'echanger la production et la consommation du trait en question. L'annulation du trait ayant \'et\'e effectu\'ee, les objets syntaxiques ne peuvent plus bouger (cf. \cite{Sta97, Sta99, Sta01}).\\
De nombreuses applications de la logique lin\'eaire et de ses syst\`emes d\'eriv\'es \`a la linguistique formelle ont d\'ej\`a eu lieu, qu'elles s'inscrivent dans le cadre des grammaires cat\'egorielles (\cite{Moo96a, Moo97, Moo99, Mor94, Oehrle95a}) ou du programme minimaliste (\cite{LecRet01, Lec03, Lec05, LecRet99, AnLec06, Amb07, Lec08}).\\
En particulier, les travaux de M. Moortgat ont mis l'accent sur une autre r\`egle sur laquelle on peut jouer dans la pr\'esentation d'un syst\`eme: la r\`egle {\it d'associativit\'e}. M\^eme si l'on juge parfois qu'un syst\`eme sans associativit\'e est de peu d'int\'er\^et math\'ematique (\cite{Gi07}), sur le plan linguistique, il fournit une sensibilit\'e \`a la {\it constituance}, autrement dit \`a la pr\'esentation possible d'une structure syntaxique sous une forme d'arbre. Que reste-t-il dans un syst\`eme logique d\'epourvu de contraction, d'affaiblissement, de permutation et d'associativit\'e? Il reste ce que Moortgat (\cite{Moo97}) appelle une {\it logique de pure r\'esiduation}, correspondant au calcul {\bf NL} de J. Lambek (\cite{Lam61}).
$$
\begin{array}{ll}
(REFL) & A\ra A \\
(TRANS) & si~A\ra B~et~B\ra C,~alors~A\ra C \\
(RES) & A\ra C/B~ssi~A\bullet B\ra C~ssi~B\ra A\lto C
\end{array}
$$
$(REFL)$ est la r\'eflexivit\'e de la relation de d\'eduction, de m\^eme que $(TRANS)$ est la transitivit\'e de cette m\^eme relation. $(REFL)$ correspond \`a ce que, dans un calcul des s\'equents, on repr\'esente par l'axiome d'identit\'e:
$$
A\vdash A
$$
et $(TRANS)$ \`a la r\`egle de coupure (ici dans sa version intuitionniste):
$$
\prooftree
\Gamma\vdash A\quad\quad \Delta, A, \Delta' \vdash B
\justifies
\Delta, \Gamma, \Delta' \vdash B
\endprooftree
$$
Si la r\`egle de coupure para{\^\i}t requise dans tout syst\`eme logique qui se respecte, pour la raison qu'int\'eress\'es au premier chef par le concept de {\it preuve}, nous ne saurions nous passer d'un moyen de les composer entre elles (ni d'un moyen pour les tranformer \'eventuellement en preuves {\it analytiques} c'est-\`a-dire sans coupures, gr\^ace au th\'eor\`eme d'\'elimination bien connu), en revanche la question peut encore se poser en ce qui concerne l'axiome d'identit\'e. C'est cette question que nous allons aborder dans cet article. Elle va nous conduire, comme nous allons le voir, vers des syst\`emes gr\^ace auxquels on peut cerner d'une mani\`ere encore plus pr\'ecise la notion concr\`ete de ressource.
\section{Un syst\`eme sans r\`egle d'identit\'e?}
\subsection{(Im)permanence de la signification}
Il faut noter au pr\'ealable que, dans tous les syst\`emes o\`u on l'utilise, l'axiome d'identit\'e va de pair avec la d\'ecomposition atomique des formules: les r\`egles concernant les symboles logiques sont telles en effet qu'il soit toujours possible de se restreindre \`a une forme de l'axiome o\`u $A$ est un atome. L'axiome d'identit\'e exprime donc une permanence de la signification des formules atomiques, et donc par extension une permanence de la signification de formules compos\'es des m\^emes atomes avec les m\^emes op\'erations logiques aux m\^emes endroits. Cette conception s'accorde bien aux math\'ematiques, pour autant qu'elles soient fond\'ees sur un tel principe d'identit\'e du symbole \`a lui-m\^eme. Elle convient peut-\^etre moins bien aux domaines extra-math\'ematiques, caract\'eris\'es par  l'emploi du langage ordinaire.
\subsection{Locativit\'e}
C'est C. Hamblin (\cite{Ham70}) qui pointe particuli\`erement bien la d\'ependance du sens par rapport \`a un {\it lieu} dans son livre sur les {\it Fallacies}. Il s'appuie sur Aristote (\cite{Aristote}) et son \'etude des sophismes qui entrent dans la classe g\'en\'erale de ceux qui r\'esultent du fait qu'une m\^eme formulation peut cacher des significations diff\'erentes. Ainsi de l'apparent syllogisme:

\it
Tous les m\'etaux sont des \'el\'ements

Le bronze est un m\'etal

Donc le bronze est un \'el\'ement

\rm
\noindent
dont le caract\`ere fallacieux est du au fait que le terme {\it m\'etal} est employ\'e dans deux sens diff\'erents, autrement dit localis\'es dans deux th\'eories diff\'erentes: la premi\`ere est la Chimie standard, la seconde une sorte de sens commun qui assimile les aliages \`a des m\'etaux.\\
De fait, un bon nombre de {\it faux arguments} proviennent du fait que des termes en apparence semblables sont en r\'ealit\'e utilis\'es dans des sens distincts, mais dire cela semble insuffisant car la notion de sens est floue. On dira plut\^ot que des termes en apparence identiques (car exprim\'es par les m\^emes suites de caract\`eres) peuvent en r\'ealit\'e relever de discours diff\'erents quand ils sont {\it positionn\'es} de mani\`eres diff\'erentes dans un discours argumentatif.\\
La logique formelle classique est mal arm\'ee pour prendre en compte de tels cas car elle est \'evidemment bas\'ee sur l'id\'ee d'identit\'e des formules: telle formule A \'enonc\'ee \`a tel moment d'un discours (ou \`a tel {\it lieu}) est identique \`a cette autre formule \'egalement d\'esign\'ee par A intervenant \`a un autre moment ou \`a un autre lieu.\\
Hamblin (p. 286) rappelle cela fort justement:
\begin{quotation}
\noindent
an approach such as that of the previous chapter [ie: the formal classical one], by locating most of the properties of the locutions in propositional letters such as 'A', 'B', 'S' and 'T', smuggles in the fiction that the question of meaning can be isolated from that of dialectical properties. When the letter 'S', say, is used twice or more in a given example it is by convention the case that it has the same meaning at each occurrence; but if meanings are to be allowed to change with context, and to be determined by the extended context, the question of whether the meaning of a given symbol changes is to be answered a posteriori and the question should not be begged by writing in an assumption of constancy.
\end{quotation}
D'o\`u l'int\'er\^et d'une conception qui dissocie les formules et les lieux.\\
Or, une logique sans r\'eflexivit\'e est une logique o\`u ce sont les lieux qui sont en premier rep\'er\'es, et les premiers \`a compter: c'est le sens qu'il faut donner \`a l'article de Girard {\it Locus Solum} (\cite{Gi01}).
\subsection{Lieux et adresses}
Soit donc un syst\`eme d'adresses.\\
On donne ici la d\'efinition de J. Y. Girard dans \cite{Gi07}:
\begin{definition}
Un {\rm biais}, notation i, j, k ..., est un entier naturel. Une {\rm ramification}, notation I, J, K, ... , est un ensemble fini de biais. Un {\rm r\'epertoire} est un ensemble quelconque de ramifications. Un {\rm locus}, ou {\rm lieu}, {\rm adresse}, notation $\sigma, \tau, \nu, \xi, ...$ est une suite finie $<i_1, i_2, ..., i_n>$ de biais. La {\rm parit\'e} d'un {\rm locus} est d\'finie comme la parit\'e de sa longueur n.
\end{definition}
\noindent
Ces adresses sont celles o\`u peuvent venir se loger des \^etres un peu particuliers, nous les appellerons formules pour commencer, m\^eme si notre ambition est de "reconstruire" les formules par la suite, en tant qu'objets jouissant de certaines bonnes propri\'et\'es de comportement vis-\`a-vis d'autres objets de la m\^eme esp\`ece (autrement dit des {\it types}, o\`u l'on retrouve ainsi l'analogie entre {\it types} et {\it formules}). Ce sont donc des {\it pr\'e-formules} ou des {\it proto-formules}. Dans l'id\'eal, on pourrait s'en passer totalement, on en a besoin ici simplement pour nous guider heuristiquement dans la conception (le {\it design}) des r\'eseaux d'adresses.
\subsection{Dessins et polarit\'e}
Un r\'eseau d'adresses est justement appel\'e un {\it dessin}. L'aide que nous recevons des {\it pr\'e-formules} qui les {\it informent} concerne la mani\`ere dont nous pouvons interpr\'eter ces r\'eseaux comme des {\it preuves}. Admettons donc que les {\it pr\'e-formules} en question sont des formules de la logique lin\'eaire, \'eventuellement enrichies d'op\'erateurs jouant sur la polarit\'e des sous-formules. On le sait depuis les travaux d'Andr\'eoli sur la recherche de preuves en logique lin\'eaire (\cite{And92}), les objets de la logique lin\'eaire sont polaris\'es. Si nous regardons en effet les r\`egles associ\'ees aux connecteurs, nous constatons que l'introduction de certains ($\wp$, $\&$) est r\'eversible \`a la diff\'erence des autres ($\ts$ et $\oplus$) dont l'introduction ne l'est pas. Les premiers sont dits {\it n\'egatifs}, les seconds {\it positifs}. Une formule dont le connecteur principal est positif (resp. n\'egatif) est dite positive (resp. n\'egative). Les r\`egles elles-m\^emes sont class\' ees en positives et n\'egatives selon qu'elles introduisent un connecteur positif ou un connecteur n\'egatif. Andr\'eoli montre qu'on peut toujours amener une preuve \`a avoir une forme canonique alternant les pas positifs et les pas n\'egatifs. 
Evidemment, cela signifie qu'\`a chaque pas polaris\'e, plusieurs introductions peuvent avoir lieu simultan\'ement, ou dit autrement, que si un connecteur est introduit, c'est un connecteur {\it synth\'etique}, la r\`egle ayant alors un nombre arbitraire de pr\'emisses ({\it logique hypers\'equentialis\'ee}).\\
De l\`a l'id\'ee qu'on peut formuler un syst\`eme tr\`es g\'en\'eral en n'utilisant que deux r\`egles "logiques": la positive et la n\'egative.\\
Ces r\`egles portent sur des objets ayant la forme g\'en\'erale de "s\'equents", mais qui ne gardent que des adresses.
\begin{definition}
Un objet $\Gamma\vdash\Theta$ o\`u $\Gamma$ est un {\rm locus} et $\Theta$ une suite de {\rm loci}, est appel\'e une {\rm fourche}, {\rm positive} si $\Gamma = \emptyset$ , {\rm n\'egative} sinon
\end{definition}

\noindent
La r\`egle positive est:
$$
\displaylines{
\prooftree
... \quad\xi\star i \vdash \Lambda_i \quad...
\justifies
\vdash \xi, \Lambda
\using
(+, \xi, I)
\endprooftree
}
$$
\begin{itemize}
\item i parcourt I
\item les $\Lambda_i$ sont {\it deux \`a deux disjoints} et inclus dans $\Lambda$
\end{itemize}
Si on lit la r\`egle du bas vers le haut, $\xi$ \'etant un \'el\'ement quelconque de la fourche, la r\`egle dit simplement que l'on {\it choisit} un lieu, avant de le distribuer sur une certaine ramification. Elle repr\'esente donc bien une action {\it positive}. Le contexte $\Lambda$ est \'eclat\'e en sous-contextes deux \`a deux disjoints, mais dont l'union ne contient pas n\'ecessairement tous les lieux du contexte d'origine.\\
La r\`egle n\'egative est:
$$
\displaylines{
\prooftree
... \quad\vdash \xi\star J, \Lambda_J \quad...
\justifies
\xi \vdash  \Lambda
\using
(-, \xi, \cal{N})
\endprooftree
}
$$
\begin{itemize}
\item J parcourt $\cal{N}$
\item les $\Lambda_J$ sont inclus dans $\Lambda$
\end{itemize}
Toujours avec la m\^eme lecture ascendante, nous constatons que le lieu utilis\'e par cette r\`egle est impos\'e: c'est le seul lieu en partie n\'egative de la fourche. Cette r\`egle repr\'esente donc une action {\it passive} ou n\'egative. Le lieu impos\'e est distribu\'e sur un r\'epertoire, et le contexte $\Lambda$ est \'eclat\'e en diff\'erent sous-contextes, qui ne forment pas n\'ecessairement une partition de $\Lambda$: ils peuvent se chevaucher, et on peut aussi perdre certains lieux en cours de route.\\
Si nous n'avons plus de liens axiomes dans les r\'eseaux, comment allons-nous les arr\^eter? Il existe une r\`egle positive particuli\`ere que Girard appelle {\it da{\"\i}mon}, sur laquelle nous reviendrons par la suite.\\
Il s'agit de la r\`egle:
$$
\displaylines{
\prooftree
\justifies
\vdash \Lambda
\using
\dagger
\endprooftree
}
$$
Elle est positive.

\noindent
Comme rappel\'e en introduction, une caract\'eristique fondamentale des syst\`emes logiques, si on veut du moins qu'ils refl\`etent des propri\'et\'es de coh\'erence du point de vue des preuves (qu'elles puissent \^etre vues comme des processus composables et gardant sous des transformations r\'egl\'ees certains invariants) est la propri\'et\'e d'\'elimination des coupures, qui exprime la r\'eelle dynamicit\'e des syst\`emes. Ici, il n'y a pas de r\`egle de coupure proprement dite, puisque nous regardons les ph\'enom\`enes au niveau des lieux et non \`a celui des formules (donc pas de formule {\it de coupure}!), mais il y a la propri\'et\'e selon laquelle une m\^eme adresse $\xi$ peut \^etre occup\'e simultan\'ement par deux instances d'un m\^eme contenu, mais avec des polarit\'es oppos\'ees. Dans ce cas, nous avons une {\it interaction} qui conduit \`a la neutralisation de cette adresse: le r\'eseau se r\'e\'ecrit au moyen du syst\`eme des sous-adresses. Si le m\^eme \'ev\`enement se reproduit au niveau des sous-adresses et ainsi de suite jusqu'\`a ce que toutes les adresses soient neutralis\'ees, on arrive n\'ecessairement sur un hypers\'equent vide, auquel on pourra alors appliquer la r\`egle du {\it da{\"\i}mon}, on dira dans ce cas que le r\'eseau a \'et\'e normalis\'e ou que sa normalisation a converg\'e. Mais pour que ce processus ait lieu, il aura fallu un r\'eseau particulier, qui ne se compose pas d'un seul {\it dessin} mais d'au moins deux {\it dessins}, dont chacun pourrait \^etre qualifi\'e de {\it contre-dessin} par rapport \`a l'autre.

\subsection{Des dessins aux desseins}
Un tel processus fait in\'evitablement penser \`a une {\it confrontation} entre deux actants\footnote{on peut aussi penser \`a l'harmonisation de deux brins d'ADN}, ou du moins deux ensembles d'actions encha{\^\i}n\'ees les unes \`a la suite des autres. \\

\noindent
Cette structure est la figure id\'eale d'un dialogue (ou, pour le dire comme Girard, d'une {\it dispute}). D'o\`u le fait que les dessins puissent aussi se voir comme des ensembles d'actions destin\'es \`a entrer en contact avec d'autres du m\^eme genre. Sous cet angle, on peut parler de {\it desseins} au lieu de {\it dessins}, et on voit surgir une interpr\'etation \`a ce qui n'\'etait jusqu'\`a maintenant que des objets syntaxiques, une interpr\'etation en termes de {\it jeux}.\\

\noindent
Dans un travail commun avec Myriam Quatrini (non publi\'e), nous avons voulu donner un exemple de cette dualit\'e preuve-strat\'egie en consid\'erant le cas de l'interpr\'etation s\'emantique \`a accorder \`a une phrase contenant plusieurs quantifieurs (ici deux). Soit la phrase:\\

(1)	{\it every linguist speaks an african language}\\

\noindent
La "signification" de (1) peut \^etre donn\'ee par la possibilit\'e d'un dialogue tel que le suivant:
\begin{enumerate}
\item celui qui soutient (1), que nous appellerons P,  se d\'eclare pr\^et \`a r\'epondre \`a toute intervention concernant un individu $d$
\item un opposant O propose un individu $f$ dont il pr\'etend qu'il est linguiste et qu'il ne conna{\^\i}t aucune langue africaine
\item P propose une langue africaine $e_f$ dont il pr\'etend que $f$ la parle
\item au m\^eme moment, O est pr\^et \`a recevoir ce genre d'affirmation
\item si O reconna{\^\i}t la validit\'e de l'affirmation de P, il peut accepter et mettre fin au dialogue
\end{enumerate}
Ce dialogue peut aussi se poursuivre plus avant:
\begin{enumerate}
\item[6.] P demande \`a v\'erifier que $f$ est bien un linguiste
\item[7.] O a la possibilit\'e d'en faire la preuve (au moyen d'une base de donn\'ees par exemple)
\item[8.] corr\'elativement, O demande \`a v\'erifier que $e_f$ est bien une langue africaine et que $f$ la parle, 
\item[9.] ce dont toujours P peut faire la preuve au moyen de {\it donn\'ees}
\end{enumerate}
Ce faisant, ce dialogue utilise des {\it faits}, autrement dit des atomes assertables au moyen d'une connaissance ext\'erieure\footnote{La ludique peut donc accueillir des donn\'ees externes et se rabattre si n\'ecessaire sur la consid\'eration d'un fait atomique comme "vrai" ou "faux"}.\\
Lorsque nous faisons cela, il semble que nous ne soyons apr\`es tout pas tr\`es \'eloign\'es de d\'emarches existant depuis de nombreuses ann\'ees, en {\it Game Theoretical Semantics} (\cite{HinKul83, Hintikka}).  Notons cependant trois diff\'erences de taille (d'autres surgiront par la suite):
\begin{itemize}
\item chaque "coup" dans le jeu est une {\it interaction} entre une {\it action positive} d'un des deux participants et une {\it action n\'egative} de l'autre. Ainsi, le pas peut \^etre franchi seulement si ce qui est avanc\'e dans l'action positive d'un des deux participants correspond aux pr\'evisions et attentes de l'autre (cette caract\'eristique est absente de la {\it GTS}),
\item il n'y a pas de r\`egle bien d\'efinie a priori associ\'ee \`a tel ou tel connecteur ou quantificateur de la logique, que suivraient les participants: nous nous contentons de suivre ce qui serait un dialogue "naturel". En particulier aucune r\`egle ne vient limiter le nombre de fois o\`u un coup pourrait \^etre rejou\'e. On peut imaginer ici que P se soit tromp\'e dans le choix d'une langue africaine parl\'ee par $f$, auquel cas il pourrait recommencer\footnote{ceci, il est vrai, n\'ecessite un dispositif un peu plus sophistiqu\'e que celui qui est ici pr\'esent\'e: il faut l'analogue des exponentielles}
\item chaque "coup" jou\'e est localis\'e en un {\it foyer} sp\'ecifique, qui est une adresse qui porte la marque de toutes les adresses par o\`u les joueurs sont ant\'erieurement pass\'es: l'historique du dialogue, depuis le d\'ebut, peut \^etre pris en compte.
\end{itemize} 
L\`a toutefois o\`u l'on s'\'ecarte le plus de la {\it GTS}, c'est dans la possibilit\'e d'autres dialogues. Il se pourrait par exemple que les deux participants ne s'entendent pas du tout sur ce qu'ils entendent par "parler une langue", ainsi le fait que {\it John parle le Ew\'e} pourrait \^etre remis en cause de bien des mani\`eres, et peut-\^etre les deux actants en seraient-ils amen\'es \`a se s\'eparer sans accord (ce qui est une autre mani\`ere de dire que leur discussion pourrait s'av\'erer interminable). Autrement dit, nous envisageons le cas o\`u des expressions comme {\it John est linguiste}, {\it la langue X est une langue africaine}, {\it John parle la langue X} ne seraient pas, en d\'epit des apparences, des atomes: leur "preuve" requerrait d'autres d\'eveloppements. Nous envisageons donc des dialogues {\it infinis}.\\
Ce qui est maintenant frappant est que les \'etapes de dialogue mentionn\'ees ci-dessus peuvent tr\`es bien s'exprimer \`a partir de la d\'ecomposition de formules de logique lin\'eaire, ou plus exactement de {\it tentatives de preuve} de ces formules. Le dialogue mentionn\'e ci-dessus pourrait \^etre associ\'e \`a la formule:
$$
S_1:\quad\&_x (\uparrow L(x)\pop \oplus_y (\uparrow A(y)\otimes \uparrow P(x,y)))
$$
et le dialogue lui-m\^eme repr\'esent\'e par les diff\'erentes phases ci-dessous:\\

\begin{tabular}{c|c}
$P$ & $O$\\
\\
\footnotesize
\shortstack{\shortstack{${\cal D}_{d'}$\\
$\vdots$\\
$\quad$}
\shortstack{${\cal D}_d$\\
$\vdots$\\
$\vdash  \downarrow L^\perp(d),\oplus_y ( \uparrow A(y)\otimes \uparrow P(d,y)))$}\hspace{1em}\shortstack{${\cal D}_{d''}$\\
$\vdots$\\
$\quad$}
 \\
$\hrulefill$\\
$(\&_x  (\uparrow L(x)\pop \oplus_y  ( \uparrow A(y)\otimes\uparrow P(x,y))))^\perp\vdash$}
&
\shortstack{
\shortstack{ 
$\downarrow L^\perp(f)\vdash$}
\hspace{1em}
\shortstack{
            $(\&_y \uparrow A(y)\pop\downarrow P^\perp(f,y))^\perp\vdash$}\\
             $\hrulefill$\\
            $\vdash \oplus_x ( \uparrow L(x)\otimes \&_y  ( \uparrow A(y)\pop  \downarrow P^\perp(x,y))$}
         \\
            \\
            
$P$ pr\^et \`a donner ds justifications &             $O$ propose un individu $f$ (affirmant \\
  pour tout individu~: $d$,$d'$,\dots&     que $f$ est un linguiste et que $f$ ne\\
 &      conna{\^\i}t aucune langue africaine)\\ \hline \\
 
$P$&$O$\\
\\

\shortstack{ 
\shortstack{$\downarrow A_{e_f}\vdash$}
\hspace{0,5em}
\shortstack{$\downarrow P_{f,e_f}\vdash$}\\
 $\hrulefill$\\
$\vdash \downarrow L^\perp(f),\oplus_y ( \uparrow A_y\otimes \uparrow P_{f,y}))$}
&
\shortstack{\shortstack{$\vdots$\\
$\vdash \downarrow A_{a'}^\perp,P_{f,a'}^\perp$}
\hspace{0,5em}
\shortstack{$\vdots$\\
$\vdash \downarrow A_{a}^\perp,P_{f,a}^\perp$}
\dots\\
$\hrulefill$\\
            $(\&_y \uparrow A_y\pop\downarrow P^\perp_{f,y})^\perp\vdash$} 
    \\
     \\
      $P$ met en \'evidence une langue $e_f$   & au m\^eme moment,, $O$ est pr\^et \\
 (affirmant que $e_f$ est une langue africaine &  \`a recevoir une telle affrirmation de  $P$,\\
  et que $f$ parle $e_f$) &  pour une certaine langue parmi  $a$, $a'$ \dots\\
            \end{tabular}

\normalsize
puis :\\
\footnotesize
\begin{tabular}{c|c}
$P$&$O$\\
\\
\shortstack{ 
\shortstack{$\downarrow A_{e_f}\vdash$}
\hspace{0,5em}
\shortstack{$\downarrow P_{f,e_f}\vdash$}\\
$\hrulefill$\\
$\vdash \oplus_y ( \uparrow A_y\otimes \uparrow P_{f,y}))$}
&
\shortstack{\shortstack{ $\hrulefill_\dag$\\
$\vdash \downarrow A_{a'}^\perp,P_{f,a'}^\perp$}
\hspace{0,5em}
\shortstack{ $\hrulefill_\dag$ \\
$\vdash \downarrow A_{a}^\perp,P_{f,a}^\perp$}
\dots\\
$\hrulefill$\\
            $(\&_y \uparrow A_y\pop\downarrow P^\perp_{f,y})^\perp\vdash$} 
     \\
\\  
\normalsize
  $P$ donne une langue africaine $a$ (ou $a'$ ou $\dots$) & \normalsize
 $O$ est pr\^et \`a accepter
            \end{tabular}\\
\normalsize           
Ici, l'occurrence de $\dag$ (le {\it da{\"\i}mon}) exprime le fait que O ne va pas plus loin dans l'\'echange.\\
           
\noindent            
Le dialogue peut se poursuivre davantage:\\

\begin{tabular}{c|c}
$P$ & $O$\\
\\
 \footnotesize{
\shortstack{\shortstack{${\cal D}_{d'}$\\
$\vdots$\\
$\quad$}
\shortstack{${\cal D}_f$\\
$\vdots$\\
$\vdash \oplus_y ( \uparrow A_y\otimes \uparrow P_{f,y}))$\\
 $\hrulefill$\\
$L_f\vdash \oplus_y ( \uparrow A_y\otimes \uparrow P_{f,y}))$\\
$\hrulefill$\\
$\vdash  \downarrow L^\perp_f,\oplus_y ( \uparrow A_y\otimes \uparrow P_{f,y}))$}
\shortstack{${\cal D}_{d}$\\
$\vdots$\\
$\quad$}\\
$\hrulefill$\\
$(\&_x  (\uparrow L(x)\pop \oplus_y  ( \uparrow A(y)\otimes\uparrow P(x,y))))^\perp\vdash$}
}
&
\footnotesize{
\shortstack{
\shortstack{$\hrulefill_\emptyset$\\
$\vdash L_f$\\
$\hrulefill$\\ 
$\downarrow L^\perp_f\vdash$}
\hspace{0,5em}
\shortstack{
            $(\&_y \uparrow A_y\pop\downarrow P^\perp_{f,y})^\perp\vdash$}\\
     $\hrulefill$ \\
            $\vdash \oplus_x ( \uparrow L_x\otimes \&_y  ( \uparrow A_y\pop  \downarrow P^\perp_{x,y})$}}\\
\\
  $P$ v\'erifie et        & $O$ peut assurer que $f$   \\
  accorde que $f$ est un  linguiste.   & est un linguiste (en tant que donn\'ee) \\

            \end{tabular}\\
            
\noindent
Nous constatons alors que les suites d'actions qui s'opposent dans le dialogue correspondent bien \`a des tentatives de fournir des preuves pour des assertions compl\'ementaires (ici: {\it tout linguiste conna{\^\i}t une langue africaine} vs {\it il existe un linguiste qui ne conna{\^\i}t aucune langue africaine}). Simplement, lorsque deux tentatives de preuve s'opposent, une seule des deux peut aboutir: elle constitue en ce cas une "vraie" {\it preuve}. Nous pourrons alors seulement consid\'erer l'autre comme une {\it contre-preuve}, c'est-\`a-dire certes une "fausse" preuve mais n\'eanmoins un objet digne d'int\'er\^et, d'o\`u le fait que nous soyons amen\'es dans le cadre de la ludique \`a nous situer dans un espace qui contient preuves et contre-preuves, plus g\'en\'eralement r\'eunies sous le nom de {\it para-preuves} ou mieux, selon l'appellation donn\'ee par P. Livet, {\it d'\'epreuves}. D'autre part, l'\'epineux probl\`eme de savoir \`a quoi peuvent bien correspondre des preuves pour des propositions atomiques est r\'esolu soit par la prise en compte de {\it donn\'ees} (cas particulier de r\`egle positive, nous y reviendrons plus loin), soit par l'admission de suites d'actions arbitrairement longues (voire infinies). Dans ce second cas, rien n'est {\it prouv\'e} \`a proprement parler: nous sommes typiquement dans la perspective des dialogues infinis tels que pouvait les concevoir N\=ag\=arjuna au 2$^{eme}$ si\`ecle de notre \`ere et qui le conduisait \`a adopter l'attitude radicalement sceptique: {\it I assert nothing} (cf. \cite{Gor08}).\\
Dans tous les cas, nous mettons en \'evidence ici la possibilit\'e de voir la signification d'une phrase comme un ensemble de {\it justifications potentielles} (par rapport \`a des demandes pr\'evisibles).
\subsection{Des desseins aux comportements}
Une troisi\`eme fa{\c c}on de repr\'esenter cette signification consiste alors \`a se d\'ebarrasser de l'\'echelon "proto-formule" fourni ici par la formule $S_1$, et \`a passer purement et simplement aux r\'eseaux d'adresses, c'est-\`a-dire aux {\it dessins}. Les sch\'emas obtenus peuvent \^etre vus comme des squelettes des preuves ci-dessus. Par exemple en ce qui concerne le point de vue de P:\\

  \begin{center}
 \footnotesize{
${\cal D}=$\shortstack{
\shortstack{${\cal D}_{d'}$\\
$\vdots$}
\hspace{1em}
\shortstack{\shortstack{ $\xi.0.3^d.7\vdash$}\hspace{2em} 
\shortstack{$\xi.0.3^d.5\vdash$}\\
$\hrulefill$\\
$ \vdash\xi.0.2^d,\xi.0.3^d$}
\hspace{1em}
\shortstack{${\cal D}_{d''}$\\
$\vdots$}\\
$\hrulefill$\\
$\xi.0\vdash$\\
$\qquad\hrulefill\qquad$\\
$\vdash\xi$}
}
\end{center}

\noindent
Les lieux ($\xi, \xi.0, \xi.0.2, \xi.0.3$ etc.) sont les localisations pr\'ecises des \'el\'ements constituants de la phrase \'etudi\'ee, d'un point de vue logique.\\
Le dessin ${\cal D}$ peut interagir avec le dessin ${\cal E}$ qui correspond \`a la contre-preuve donn\'ee par O.
$$
  \displaylines{
  \footnotesize
   {\cal E} =
 \prooftree
\[
\[
\justifies
\xi.0.2^d\vdash
\using
\emptyset
\]
\[
\[
\justifies
\vdash\xi.0.3^d.7, \xi.0.3^d.5
\using
\dag
\]
\justifies
\xi.0.3^d\vdash
\]
\justifies
\vdash\xi.0
\]
\justifies
\xi\vdash
\endprooftree
  }
$$
Le $d$ en exposant indique le choix d'un individu au pas consid\'er\'e. $\emptyset$ indique que O pose que "$d$ est un linguiste" comme un fait ou une donn\'ee: P n'a d\'esormais plus la possibilit\'e de jouer sur le locus $\xi.0.2$ s'il reconna{\^\i}t ce fait.\\
L'interaction entre ${\cal D}$ et ${\cal E}$ (co{\"\i}ncidence des lieux de polarit\'es oppos\'ees dans les deux dessins) conduit \`a une neutralisation desdits lieux et au dessin \'el\'ementaire:
$$
\prooftree
\justifies
\quad\quad\vdash\quad\quad
\using
\dag
\endprooftree
$$
Ce dessin \'el\'ementaire peut \^etre per{\c c}u aussi bien comme une {\it preuve} que comme une {\it contre-preuve}, c'est une sorte de vecteur nul de l'espace des \'epreuves. Deux dessins dont l'interaction conduit \`a lui sont donc dits {\it orthogonaux}.
$$
{\cal D}\perp{\cal E}
$$
Nous noterons $[[~, ~]]$ cette sorte de {\it produit scalaire} entre dessins qui appara{\^i}t lorsque nous les faisons interagir (on parle aussi de {\it normalisation} du r\'eseau qu'ils forment).\\
Dans le cas d'orthogonalit\'e, nous avons:
$$
[[{\cal D}, {\cal E}]] = \prooftree \justifies \quad\vdash\quad \using\dag \endprooftree
$$
On peut \'ecrire aussi que: ${\cal E}\in{\cal D}^\perp$. On a bien \'evidemment: ${\cal D}\in{\cal D}^{\perp\perp}$, mais ce dernier ensemble (le bi-orthogonal de ${\cal D}$) contient bien d'autres dessins : on dira qu'il est engendr\'e par un dessin (ou un ensemble de dessins), ici ${\cal D}$.\\
On appelle {\it comportement engendr\'e par} ${\cal D}$ l'ensemble ${\cal D}^{\perp\perp}$. \\
Autrement dit, un {\it comportement} est un ensemble de dessins qui se comportent tous de la m\^eme mani\`ere vis-\`a-vis de l'interaction avec d'autres dessins. {\bf C'est par la notion de comportement que nous allons regagner la v\'eritable notion de formule}.\\
Par exemple, l'\'enonc\'e {\it John parle le Ew\'e}, s'il est associ\'e \`a au moins un dessin, engendrera un comportement ${\cal C}$ qui contiendra des \'el\'ements de plus en plus sp\'ecifi\'es au fur et \`a mesure que nous le ferons interagir avec d'autres dessins. De cette fa{\c c}on, nous pouvons associer \`a chaque \'enonc\'e \'el\'ementaire un comportement.\\
Au stade de reconstruction de la logique, nous pourrons d\'efinir des {\it op\'erations} entre comportements parmi lesquelles nous retrouverons les connecteurs de la logique lin\'eaire, de sorte que la {\it proto-formule} $S_1$ correspondra bien finalement \`a une formule, mais au lieu de la d\'efinir plus ou moins arbitrairement comme nous l'avons fait, cette formule d\'ecoulera des op\'erations {\it g\'eom\'etriques} associ\'ees aux dessins vus comme des suites d'actions. Ce sera, en un sens, une composition de dialogues \'el\'ementaires.\\

\noindent
Nous pouvons donc d\'esormais consid\'erer nos "pr\'e-formules" comme des ensembles de demandes de justifications et de justifications \'el\'ementaires. Ces ensembles, ou {\it dessins} peuvent \^etre, \`a vrai dire, arbitrairement enrichis: il en r\'esulte des ensembles de dessins dont la coh\'erence interne est simplement exprim\'ee par le fait qu'ils interagissent de la m\^eme mani\`ere avec d'autres dessins. Par cl\^oture, on obtient des comportements, qui correspondent aux vraies formules.
\subsection{La d\'elocalisation par ${\cal F}ax$}
Il reste encore, bien entendu, \`a \'etablir une possibilit\'e de {\it transfert}: la signification d'un \'enonc\'e ou d'un terme ne saurait \^etre purement li\'ee \`a des lieux, car cela voudrait dire que nous avons sans cesse \`a r\'e-inventer la langue et ses significations. Si les significations sont surtout "dialectiques" et d\'ependantes des contextes, il n'en reste pas moins que la plupart du temps, les locuteurs s'entendent sur le sens des mots, ce qui ne serait gu\`ere possible s'ils devaient le r\'einventer tout au long des usages. Au cours de l'encha{\^\i}nement des discours se cr\'eent des r\'egularit\'es r\'eutilisables dans d'autres contextes. Hamblin parle \`a ce propos de {\it patterns of use}. Le fait qu'une certaine permanence de signification puisse exister pour certains termes n'est pas incompatible avec l'id\'ee du caract\`ere profond\'ement dialectique de la signification. \\ 
Comme le dit encore Hamblin (p. 295):
\begin{quotation}
\noindent
we may have to say that in so far as there is a presumption that W is constant in meaning there is a presumption that any given use of W is part of a pattern, or that the user's explanations of his meaning are mutually coherent.
\end{quotation}

\noindent
Si nous n'avons plus d'axiome d'identit\'e (comme nous n'avons plus d'ailleurs de r\`egle de coupure explicite), comment allons-nous v\'erifier que les contenus figurant \`a des adresses diff\'erentes sont en r\'ealit\'e identiques, un contenu ayant seulement subi un transfert d'un lieu vers un autre?\\
C'est ici qu'intervient un dessin particulier: le ${\cal F}ax$. C'est un dessin infini r\'ecursivement d\'efini par:
$$
\displaylines{
Fax_{\xi, \xi'} =
\prooftree
...
\[
...
\[Fax_{\xi'_{i}, \xi_{i}}
\justifies
\xi'\star i\vdash \xi\star i
\]
...
\justifies
\vdash \xi\star J, \xi' 
\using
(+, \xi', J)\]
...
\justifies
\xi \vdash \xi'
\using
(-, \xi, {\cal P}_f({\mathbb N}))
\endprooftree
}
$$
Au premier pas, qui est n\'egatif, le lieu n\'egatif est distribu\'e sur tous les sous-ensembles finis de ${\mathbb N}$, puis, pour chaque ensemble d'adresses (relatif \`a un $J$), le lieu positif $\xi'$ est choisi et il se cr\'ee une sous-adresse $\xi'\star i$ pour chaque $i\in J$, et le m\^eme m\'ecanisme est relanc\'e pour le nouveau lieu ainsi obtenu.\\
On peut alors voir le r\^ole que joue ce dessin dans sa normalisation avec un dessin arbitraire.
Prenons un dessin ${\cal D}$ de base $\vdash\xi$, on peut facilement montrer que sa normalisation avec ${\cal F}ax_{\xi\vdash\rho}$ a pour r\'esultat ${\cal D}'$ qui n'est autre que ${\cal D}$ mais o\`u l'addresse $\xi$ est syst\'ematiquement remplac\'ee par $\rho$.\\
A titre d'exemple, prenons pour ${\cal D}$ le dessin :
$$
\displaylines{
\prooftree
\[
{\cal D}_1
\justifies
\xi\star 1\vdash
\]
\[
{\cal D}_2
\justifies
\xi\star 2\vdash
\]
\justifies
\vdash\xi
\endprooftree
}
$$
La normalization a lieu en s\'electionnant la "tranche" correspondant au sous ensemble $\{1, 2\}$ de sorte qu'une fois \'elimin\'ee la premi\`ere coupure, il reste :
$$
\displaylines{
\hfill
\prooftree
{\cal D}_1
\justifies
\xi\star 1\vdash
\endprooftree
\hfill
\prooftree
{\cal D}_2
\justifies
\xi\star 2\vdash
\endprooftree
\hfill
\prooftree
\[
{\cal F}ax
\justifies
\rho\star 1\vdash\xi\star 1
\]
\quad
\[
{\cal F}ax
\justifies
\rho\star 2\vdash\xi\star 2
\]
\justifies
\vdash\xi\star 1, \xi\star 2, \rho
\endprooftree
\hfill
}
$$
Les deux dessins de gauche normalisent avec celui de droite, donnant finalement:
$$
\displaylines{
\prooftree
\[
{\cal D}_1'
\justifies
\rho\star 1\vdash
\]
\[
{\cal D}_2'
\justifies
\rho\star 2\vdash
\]
\justifies
\vdash\rho
\endprooftree
}
$$
o\`u, dans ${\cal D}_1'$ and ${\cal D}_2'$, l'addresse $\xi$ est syst\'ematiquement remplac\'ee par $\rho$.\\
Nous en concluons que la normalisation avec ${\cal F}ax$ est la proc\'edure qui d\'eplace un dessin d'une localisation donn\'ee \`a une autre: c'est une proc\'edure de {\it d\'elocalisation}.
\subsection{De quelques comportements \'el\'ementaires}
Parmi tous les desseins concevables, il en est certains qui sont remarquables notamment par leur bri\`evet\'e. Nous avons d\'ej\`a fait r\'ef\'erence, dans l'exemple pr\'ec\'edent, au cas o\`u l'un des partenaires emp\^eche l'autre de poursuivre autrement qu'en renon{\c c}ant \`a la continuation de la dispute ou en entrant clairement dans un cas de divergence (\emph{dissensus}). Cela peut arriver aussi bien au cours de l'accomplissement d'un pas n\'egatif que d'un pas positif. Dans le premier cas, cela revient \`a prendre $\cal{N}=\emptyset$ dans la r\`egle n\'egative. Le partenaire ne peut alors effectuer aucun mouvement sans sombrer dans le dissensus (puisqu'aucune ramification n'est jouable). Girard appelle {\it sconse} ce dessein, \`a cause de son caract\`ere particuli\`erement "associal" (!).
$$
\displaylines{
\prooftree
\justifies
\xi\vdash
\using
(\xi, \emptyset)
\endprooftree
}
$$
Dans le cas positif, nous pouvons consid\'erer le dessin suivant, qui est positif:
$$
\displaylines{
\prooftree
\justifies
\vdash\xi
\using
(\xi, \emptyset)
\endprooftree
}
$$
Si le partenaire n\'egatif vise \`a un consensus, il ne peut r\'epondre que par le \emph{da{\"\i}mon} n\'egatif puisque le proposant ne lui donne aucune adresse pour encha{\^\i}ner. Il s'agit de:

$$
\displaylines{
\prooftree
\[\justifies
\vdash
\using
\dagger
\]
\justifies
\xi\vdash
\endprooftree
}
$$
Ce faisant, il s'avoue vaincu (du moins, dans la majeure partie des interpr\'etations dont le \emph{da{\"\i}mon} peut se rev\^etir). Ainsi, \`a ce jeu, le partenaire positif gagne toujours. Evidemment, dans une situation r\'eelle, l'usage de ce dessein ne sera autoris\'e que sous des conditions r\'eguli\`eres (par exemple l'existence d'un fait comme donn\'ee). Girard nomme \emph{Bombe atomique} cet argument dissuasif! Nous le noterons aussi ${\cal B}ombe^+$. Son unique dessein orthogonal est donc le ${\cal D}${\it a{\"\i}mon} n\'egatif donn\'e ci-dessus. Il est \'evident que celui-ci est aussi orthogonal \`a:
$$
\displaylines{
\prooftree
\justifies
\vdash\xi
\using
\dagger
\endprooftree
}
$$
\noindent
Nous obtenons donc un comportement qui contient {\it deux} desseins: ${\cal B}omb^+$ et {\it da{\"\i}mon}. Notons-le {\bf 1}: il deviendra l'\'el\'ement neutre de $\ts$, tel que red\'efini en termes ludiques.\\
Consid\'erons maintenant de nouveau le {\it sconse}. Notons $\top$ le comportement n\'egatif qui le contient. Quel est l'orthogonal de $\top$? La seule possibilit\'e pour l'adversaire est de jouer le {\it ${\cal D}$a{\" \i}mon}, positif, cette fois. D'o\`u:
$$
\top^\bot = \{Dai\}
$$
et
$$
\top = Dai^\bot
$$
autrement dit, {\it tous} les desseins n\'egatifs de m\^eme base.
\subsection{Le tour de force de la ludique: regagner les op\'erateurs logiques}
Nous n'entrerons pas ici dans le d\'etail de la reconstruction de la logique op\'er\'ee par Girard \`a partir des op\'erations d\'efinissables sur les comportements. Disons simplement que, de m\^eme que de telles op\'erations sont d\'efinissables sur les {\it espaces coh\'erents}, elles le sont sur les comportements. Il est ainsi possible de d\'efinir une op\'eration, not\'ee $\odot$, entre deux desseins, puis  \`a partir de l\`a, l'op\'eration $\ts$ entre deux comportements.
\begin{definition}
Soit $\cal{U}$ et $\cal{B}$ deux desseins positifs, nous d\'efinissons le produit tensoriel $\cal{U}\odot\cal{B}$ par:
\begin{itemize}
\item si l'un des deux est le {\it ${\cal D}$a{\" \i}mon}, alors ${\cal{U}}\odot{\cal{B}} = Dai$,
\item sinon, soit $(+, <>, I)$ et $(+, <>, J)$ les premi\`eres actions de respectivement $\cal{B}$ et $\cal{U}$, si$I\cap J\not=\emptyset$, alors ${\cal{U}}\odot{\cal{B}}$ = $Dai$. Sinon, on remplace dans chque chronique de $\cal{B}$ et de $\cal{U}$ la premi\`ere action par $(+, <>, I\cup J)$, donnant respectivement $\cal{B'}$ and $\cal{U'}$, alors $\cal{U}\odot\cal{B}=\cal{U'}\cup\cal{B'}$.
\end{itemize}
\end{definition}
Il est alors possible de d\'efinir le $\otimes$  de deux comportements au moyen de d\'elocalisations. Notons toutefois que la Ludique  permet de d\'efinir de nouveaux connecteurs. \\
On peut d\'efinir le produit  $\odot$ de deux comportements par : ${\bf F} \odot {\bf G} = \{{\cal A}\odot{\cal B}; {\cal A}\in{\bf F}, {\cal B}\in{\bf G}\}$ (cela ne donne pas n\'ecessairement un comportement), et alors le produit tensoriel proprement dit par ${\bf F} \otimes {\bf G} = ({\bf F} \odot {\bf G})^{\bot\bot}$.\\
Si nous revenons alors \`a l'exemple du paragraphe 2.5, il appara{\^\i}t que le comportement associ\'e \`a $S_1$ peut s'\'ecrire:
\begin{center}
${\mathbb S}_1=\&_x (\downarrow {\mathbb L}(x)\linefle \oplus_y(\downarrow {\mathbb A}(y)\otimes \downarrow {\mathbb P}(x,y)))$
\end{center} 
construit \`a partir des comportements $ {\mathbb L}(x)$ , ${\mathbb A}(y)$ et $ {\mathbb P}(x,y)$!
Ces comportements \'el\'emen\-taires peuvent \^etre simplement remplac\'es par {\bf 1} dans une interpr\'etation factuelle de ce genre d'\'enonc\'e, mais ils peuvent aussi \^etre remplac\'es par des comportements beaucoup plus riches si le dialogue est destin\'e \`a durer, voire des comportements contenant des dessins infinis.
\section{Pragmatique et sophismes}
\subsection{Actes de langage}
Selon la th\'eorie des actes de langage (\cite{Austin, Searle}), certaines \'enonciations ({\it questions, ordres, promesses, actes de juridiction etc.}) se distinguent particuli\`erement par les transformations qu'elles font subir au contexte. Par exemple, une question provoque, dans les cas de consensus minimal, une obligation de r\'eponse, un ordre provoque, dans les m\^emes conditions, une obligation d'acqui\`escement, une promesse un changement dans les engagements auxquels sont soumis les locuteurs etc. La ludique se pr\^ete \`a leur repr\'esentation dans la mesure o\`u 
l'utilisation du ${\cal F}ax$ dans un dessin est un moyen de reproduire de telles transformations du contexte. Notons ici que nous pouvons entendre par {\it contexte} un r\'eseau de lieux. Nous pouvons dire que, dans un tel r\'eseau, certains sous-r\'eseaux sont {\it activ\'es} en fonction de ceux qui sont choisis pour loger des contenus particuliers. Nous risquerons ici une m\'etaphore neuronale: de m\^eme que les neurones individuels sont suppos\'es \^etre des adresses arbitraires pour des contenus, l'activit\'e neuronale proprement dite se rep\'erant aux sous-r\'eseaux activ\'es, on peut consid\'erer un contexte comme un ensemble arbitraire de {\it loci} et assimiler un ensemble activ\'e d'adresses \`a un {\it \'etat mental}. En ce cas, nous d\'ecrivons un acte de langage comme une transformation d'\'etat mental.\\
L'exemple le plus simple concerne le jeu {\it question-r\'eponse}. Reprenons l'exemple du paragraphe 2.5, et opposons \`a l'assertion $S$ = {\it tout linguiste conna{\^\i}t une langue africaine}, la question $Q$ = {\it Quelle langue africaine parle John?}.

\begin{tabular}{c|c}
 \footnotesize{
\shortstack{\shortstack{${\cal D}_{d'}$\\
$\vdots$\\
$\quad$}
\shortstack{\shortstack{ 
\shortstack{$\hrulefill_\emptyset$\\
$\vdash A_{e_{j}}$\\
$\hrulefill$\\
$ A _{e_{j}}^\bot\vdash$}
\hspace{0,5em}
\shortstack{$\hrulefill_\emptyset$\\
$\vdash S_{john,e}$\\
$\hrulefill$\\
$S_{john,e_{j}}^\bot\vdash$}\\
$\hrulefill$\\
$\vdash \oplus_y ( A_y\ts S_{john,y}))$}\\
$\hrulefill$\\
$L_{john}\vdash \oplus_y (  A_y\ts   S_{john,y}))$\\
$\hrulefill$\\
$\vdash  L_{john}^\bot,\oplus_y (  A_y\ts  S_{john,y}))$}
\shortstack{${\cal D}_{d''}$\\
$\vdots$\\
$\quad$}\\
$\hrulefill$\\
$M_1^\bot\vdash$}
}
&
\footnotesize{
\shortstack{
\shortstack{$\hrulefill_{\emptyset}$\\
$\vdash L_{john}$\\
$\hrulefill$\\
$L_{john}^\bot\vdash$}\hspace{0,5em}\shortstack{ \dots\hspace{0,2em}
\shortstack{ ${\cal F}ax$\\
$A_{e_j}\vdash  A_e$\\
$\hrulefill$\\
$\vdash  A_e, A_{e_j}^\bot$\\
$\hrulefill$\\
$S_{john,e}\vdash A_e, A_{e_j}^\bot$\\
$\hrulefill$\\
$\vdash  A_e, A_{e_c}^\bot,S_{john,e}^\bot$\\
$\hrulefill$\\
$ A_e^\bot\vdash A_{e_c}^\bot, S_{john,e_j}^\bot$\\
$\hrulefill$\\
$\vdash A_{e_j}^\bot,S_{john,e_j}^\bot  , A$}\hspace{1em}\dots\hspace{1em}\\
$\hrulefill$\\
$(\&_y ( A(y) \linefle S_{john,y}^\bot)^\bot\vdash A$}\\
$\hrulefill$\\
$\vdash M_1^\bot, A$}
}

            \end{tabular}\\
            
\noindent           
La proto-formule $A$ est introduite dans le contexte comme {\it r\'eceptrice} du r\'esultat de l'interaction. Apr\`es normalisation du r\'eseau, il reste seulement:

\begin{center}
\footnotesize{
\shortstack{$\hrulefill_\emptyset$\\
$\vdash A_{e_j}$\\
$\hrulefill$\\
$\vdash A=A_1\oplus\dots\oplus A_{e_j}\oplus\dots$}
}
\end{center}
o\`u  $A_{e_j}$ est la langue africaine que parle John.\\
D'autres exemples sont fournis dans un article non encore publi\'e de M-R. Fleury et S. Tron{\c c}on (\cite{FleuTron}).

\subsection{Du traitement des sophismes}
De ce qui pr\'ec\`ede, ressort l'id\'ee que si, en nous passant de r\`egles structurelles comme la contraction et l'affaiblissement, nous gagnons une sensibilit\'e aux ressources du point de vue de leur quantit\'e disponible, si en nous privant de la commutativit\'e, nous gagnons une sensibilit\'e \`a l'ordre dans lequel elles sont pr\'esent\'ees (et en nous privant de l'associativit\'e, \`a la structure en constituants), en supprimant l'axiome d'identit\'e, nous avons acquis plus encore: une sensibilit\'e aux adresses ou aux lieux. D\'esormais, nous pouvons jouer sur ces lieux comme nous avons jou\'e pr\'ec\'edemment sur le caract\`ere transitoire ou p\'erenne d'une ressource. C'est ce point qui va nous appara{\^\i}tre comme pr\'ecieux dans la discussion sur les sophismes.\\
L'un des sophismes les plus connus l'est sous le nom de {\it p\'etition de principe}. Il consiste, selon Aristote dans les {\it R\'efutations sophistiques} \`a {\it faire entrer en ligne de compte dans les pr\'emisses la proposition initiale \`a prouver} (Organon, VI, 6) (\cite{Aristote}) autrement dit \`a prouver une th\`ese qui figure elle-m\^eme d\'ej\`a, implicitement, dans les pr\'emisses du raisonnement. Mais, comme beaucoup d'auteurs l'ont not\'e (\`a commencer par l'illustre John Maynard Keynes), n'est-ce pas l\`a justement ce qu'on fait toujours en logique formelle? Autrement dit une des raisons pour lesquelles, comme l'affirmait Wittgenstein, la logique ne serait concern\'ee que par des tautologies, c'est-\`a-dire des v\'erit\'es sans contenu, valides seulement par leur forme? De fait, la logique classique (aussi bien qu'intuitionniste etc.) nous est de peu d'utilit\'e pour \'etudier cette figure de raisonnement, l'axiome d'identit\'e $A\vdash A$ semblant \^etre la forme la plus raccourcie et condens\'ee de la p\'etition de principe! \\
Si l'on veut proposer une analyse de la p\'etition de principe, il semble donc qu'un syst\`eme ne poss\'edant pas l'axiome d'identit\'e soit requis. Une phrase comme {\it l'\^ame est immortelle parce qu'elle ne meurt jamais} (\cite{Qua08}) peut \^etre d\'ecrite \`a partir d'une localisation en un lieu $\xi$ de la phrase {\it l'\^ame est immortelle} et d'une unique justification: ({\it l'\^ame ne meurt jamais}). Lorsque le locuteur P se place en $\xi$, pour d\'emarrer le dessin ${\cal D}_\xi$ de son argumentaire, il ex\'ecute, par son \'enonciation, une action {\it positive} qu'il s'appr\^ete \`a justifier. Autrement dit, tout de suite apr\`es, il commet une action n\'egative par laquelle il se montre pr\^et \`a recevoir une demande de justification (en $\xi.1$). Il poss\`ede, d\'ej\`a toute pr\^ete, cette justification: elle tiendrait en une seconde \'enonciation, {\it l'\^ame ne meurt jamais}, elle-m\^eme suppos\'ee poss\'eder sa propre justification. En affirmant ce deuxi\`eme \'enonc\'e, P d\'emarre en fait un deuxi\`eme dessin, au lieu $\xi.1.1$ (${\cal D}_{\xi.1.1}$), mais parce que le deuxi\`eme \'enonc\'e est d\'ej\`a contenu dans le premier, ${\cal D}_{\xi.1.1}$ n'est autre que le d\'ecalage par ${\cal F}ax$ du dessin ${\cal D}_\xi$! Autrement dit:
$$
{\cal D}_{\xi.1.1} = [[{\cal D}_\xi, {\cal F}ax_{\xi, \xi.1.1}]]
$$
De m\^eme que, d'ailleurs:
$$
{\cal D}_\xi = [[{\cal D}_{\xi.1.1}, {\cal F}ax_{\xi.1.1, \xi}]]
$$
Or, une th\'eorie de l'argumentation interdirait ce genre de circularit\'e.\\
Nous voyons ici que, de m\^eme que l'introduction de l'exponentielle en logique lin\'eaire "explique", par le passage qu'elle op\`ere du parfait \`a l'imparfait, la production de figures paradoxales, l'introduction du ${\cal F}ax$ permet de d\'ecrire le fonctionnement d'une autre figure, qui est une sorte {\it d'envers} du paradoxe puisqu'il s'agit au contraire d'"\'evidences" trop triviales: la p\'etition de principe. \\
De m\^eme que sans exponentielle, le {\it Menteur} se r\'esoudrait en une suite infinie d'instants altern\'es o\`u sont vrais tour \`a tour: {\it je mens} et {\it je dis la v\'erit\'e}, sans ${\cal F}ax$, la p\'etition de principe se r\'esoudrait \`a une reproduction ind\'efinie du m\^eme.\\


\noindent
Un autre sophisme fameux est celui qu'Aristote et sa post\'erit\'e ont qualifi\'e de 
sophisme de {\it plusieurs questions}. Aristote (\cite{Aristote}) y r\'ef\`ere comme au cas o\`u il y a une pluralit\'e de questions qui demeure inaper{\c c}ue et {\it qu'on donne une seule r\'eponse comme s'il n'y avait qu'une seule question}. Ce sophisme est parfois illustr\'e par le dialogue suivant:

\it
- avez-vous cess\'e de battre votre p\`ere?

- oui (ou : non)

\rm
\noindent
Dans les deux cas, le r\'epondant est reconnu coupable de brutalit\'e envers son p\`ere, que celle-ci soit pass\'ee ou pr\'esente. Le dessin de P (o\`u l'on voit que l'orthographe {\it dessein} est parfois la plus appropri\'ee!) contient alors un lieu escamot\'e, sur lequel O ne peut pas r\'epondre, d'o\`u, pour que le processus de normalisation converge, la seule possibilit\'e qu'il r\'eponde "oui" ou "non" \`a la question principale.

\begin{center}
 \shortstack{
 \shortstack{ 
 	    ${\cal L}oc_1$~: \\ ${\cal D}$\\
	    $\hrulefill$\\
             $B^\bot\vdash$\\
                       $\hrulefill_W$\\
                       $\vdash A^\bot, B$\\
                       $\hrulefill$\\
             $A\ts B^\bot\vdash$\\
             $\hrulefill$\\
             $\vdash A~\pop B$}
             \hspace{2em}
             \shortstack{
             ${\cal L}oc_2$~:\\
             $\hrulefill$\\
             $[\vdash A]\quad\quad$
             $\vdash B^\bot$\\
             $\hrulefill$\\
             $A~\pop B\vdash$}\\
             $\hrulefill$}
             \end{center}
             
\noindent
Nous nous sommes aid\'es de proto-formules pour construire ce r\'eseau de deux dessins. La formule mise au foyer $\xi$ de ${\cal L}oc_1$ est une implication lin\'eaire: normalement, pour acc\'eder au lieu de la question principale ({\it avez-vous cess\'e de battre votre p\`ere?}) il faut passer obligatoirement par la sous-question {\it battiez-vous votre p\`ere auparavant?}, autrement dit, il faudrait "consommer" la premi\`ere question avant de r\'epondre \`a la deuxi\`eme, d'o\`u $A\pop B$. Mais le locuteur \'enon{\c c}ant cette question ayant pour dessein d'escamoter $A$, se pr\'epare \`a "justifier" sa question au moyen d'un dessin qui comporte un affaiblissement sur $A$ (pas marqu\'e par $W$). Le lieu de $A$ ayant disparu, toute intervention de ${\cal L}oc_2$ sur ce lieu fera diverger l'interaction: c'est bien ce qui se passe lorsque, dans un dialogue, l'un des participants remet en cause un pr\'esuppos\'e.\\

\noindent
On peut maintenant remarquer que, si nous modifions l\'eg\`erement le jeu, avec cette fois, une formule quelconque $A\wp B_1 \wp ... \wp B_n$ au foyer du premier locuteur, o\`u chaque sous-formule peut a priori \^etre remise en cause, si ${\cal L}oc_1$  sait que ${\cal L}oc_2$ poss\`ede une strat\'egie gagnante (par exemple ${\cal B}ombe^+$) sur l'une des sous-formules, ici $A$, il ne peut, s'il veut poursuivre le dialogue, que jouer sur les $B_i$. Ceci est une sorte d'envers de la pr\'esupposition : ${\cal L}oc_1$ se prive lui-m\^eme d'un lieu dans un dessein par anticipation sur la r\'eponse de ${\cal L}oc_2$.
\begin{center}
 \shortstack{
 \shortstack{ 
 	    $(B_1~\wp~...\wp~ B_n)^\bot\vdash$\\
             $\hrulefill$\\
             $\vdash A~\wp~ B_1~\wp~...\wp~ B_n$}
             \hspace{2em}
             \shortstack{
             $\hrulefill_\emptyset$\quad\quad\quad\quad\quad\quad\quad\quad\quad\quad\\
             $\vdash A\quad\quad$
             $\vdash B_1\quad ... \quad \vdash B_n$\\
             $\hrulefill$\\
             $A~\wp~ B_1~\wp~...\wp~ B_n\vdash$}\\
             $\hrulefill$}
             \end{center}
\section{Conclusion}
Nous avons montr\'e dans cet article qu'il \'etait possible de jouer sur l'axiome d'identit\'e (la r\'eflexivit\'e de la relation de d\'eduction) comme il est possible de jouer sur les r\`egles structurelles d'un syst\`eme. Cela ne signifie pas seulement la suppression de ces r\`egles, mais leur remplacement, la plupart du temps, par un dispositif plus souple qui permet l'analyse fine de ph\'enom\`enes: paradoxes dans les cas des r\`egles structurelles, p\'etition de principe dans le cas de l'axiome d'identit\'e. Une cons\'equence fondamentale de cette \'elimination de l'identit\'e est la remise en cause de la notion de formule en tant qu'\^etre id\'eal et "spiritualiste" (au sens de Girard), au profit des \emph{adresses} o\`u peuvent se loger des \emph{contenus}. Avoir certaines proc\'edures comme le ${\cal F}ax$ permet d'assurer la d\'elocalisation (ou le transfert) de ces contenus. La notion de formule (et donc de \emph{sens} d'une formule) est regagn\'ee via la notion d'interaction, qui n'est rien d'autre que \emph{l'usage}: un \emph{dessin} s'utilise en le faisant interagir avec d'autres. L'interpr\'etation des dessins en termes de strat\'egies (\emph{desseins}) ouvre la voie \`a un approfondissement de la notion de \emph{jeu de langage} au sens de Wittgenstein, en \'evitant les inconv\'enients de la \emph{Game Theoretical Semantics}, dans la mesure o\`u la notion de jeu qui \'emerge est beaucoup plus g\'en\'erale que celle propos\'ee par Hintikka et al. : elle ne suppose ni strat\'egie gagnante, ni fonction de gain.

\bibliography{lecomte}
\bibliographystyle{plain}
\end{document}